\documentclass[a4paper, 10pt]{amsproc}
\usepackage{defs}
\usepackage{amsaddr}
\usetikzlibrary{matrix}
\usepackage{derivative}
\usepackage[pagewise]{lineno}
\usepackage{orcidlink}
\usepackage{etoolbox}

\title[DatDri FeedLearn]{Data-driven learning of feedback maps for explicit robust predictive control: an approximation theoretic view}

\author[S. Ganguly, S. Gupta, D. Chatterjee]{Siddhartha Ganguly\,\orcidlink{0000-0003-2046-2061}, \,Shubham Gupta,\,Debasish Chatterjee\,\orcidlink{0000-0002-1718-653X}}
\thanks{%
	S. Ganguly is with \faGroup\ The Applied Mathematics and Physics Department, Graduate School of Informatics, \faUniversity\ Kyoto University, \faMapMarker\  Kyoto, Japan. S. Gupta is with \faGroup\ Department of Mechanical and Process Engineering \faUniversity\ ETH Zurich, \faMapMarker\ Switzerland. D. Chatterjee is with \faGroup\ Systems \& Control Engineering, \faUniversity\ IIT Bombay, Powai, \faMapMarker\ Mumbai 400076, India.}

\thanks{%
	Contact Information: (SG) \faHome\ \url{https://sites.google.com/view/siddhartha-ganguly}, \faEnvelope\ \texttt{ganguly.siddhartha.7p@kyoto-u.ac.jp}. (SuG) \faEnvelope\ \texttt{shugupta@ethz.ch}. (DC) \faHome\ \url{https://www.sc.iitb.ac.in/~chatterjee}, \faEnvelope\ \texttt{dchatter@iitb.ac.in}
}

\date{\DTMnow}

\begin{document}
\subjclass[2020]{49N35, 93B51, 93B52, 90C34, 62M45, 41A05} 
\keywords{Robust optimal control, feedback synthesis, model predictive control, approximation theory, neural networks.}

\maketitle

\begin{abstract}
We establish an algorithm to learn feedback maps from data for a class of robust model predictive control (MPC) problems. The algorithm accounts for the approximation errors due to the learning directly at the synthesis stage, ensuring recursive feasibility by construction. The optimal control problem consists of a linear noisy dynamical system, a quadratic stage and quadratic terminal costs as the objective, and convex constraints on the state, control, and disturbance sequences; the control minimizes and the disturbance maximizes the objective. We proceed via two steps --- \textbf{(a)} \emph{Data generation:} First, we reformulate the given minmax problem into a convex semi-infinite program and employ recently developed tools to solve it in an \emph{exact} fashion on grid points of the state space to generate \texttt{(state, action)} data. \textbf{(b)} \emph{Learning approximate feedback maps:} We employ a couple of approximation schemes that furnish tight approximations within preassigned \emph{uniform} error bounds on the admissible state space to learn the unknown feedback policy. The stability of the closed-loop system under the approximate feedback policies is also guaranteed under a standard set of hypotheses. Two benchmark numerical examples are provided to illustrate the results. 

\end{abstract}




\newcommand{\paramh}{h}
\newcommand{\paramD}{\mathcal{D}}

\section{Introduction}
\label{s:intro}

This article tackles the problem of learning constrained feedback policies for robust model predictive control (MPC) of uncertain linear systems using offline data. Specifically, it explores how offline data, which may include historical measurements or pre-collected input-output trajectories, can be leveraged to train a function approximator that generates (approximate) feedback policies for constrained robust optimal control. Once trained offline, the feedback policy is implemented online, where it continuously responds to the system's states in order to undertake control decisions under constraints in real time. To this end, we focus on the task of performing receding horizon control (RHC), commonly known as model predictive control, in a fast and efficient manner by learning feedback maps/policies from offline data (hence, data-driven) with uniform approximation guarantees.

MPC \cite{ref:XiLi-19} is a dynamic optimization/optimal control technique for synthesizing constrained feedback control that has been successfully applied across an array of industries including electrical \cite{ref:MPC_for_energy}, robotics and trajectory generation \cite{ref:MPC-AUV-kamel2017model, ref:SGRADDC-24}, chemical \cite{ref:MPC-chem-marquez2019model}, power electronics \cite{ref:MPC_for_power_elec}, oil and gas, among several others. The ability to directly integrate constraints into the control design process sets MPC apart from other conventional control techniques; consequently, aggressive control actions are generated and the closed-loop processes tend to operate at the boundary of their admissible regions --- an extremely desirable feature in applications.


Let us take a high-level look at the MPC framework: In MPC, the ``model'' represents a discrete-time controlled dynamical system defined on a Euclidean space, formulated as \(x_{t+1}= f(x_t, u_t)\) for all \(t\in \Nz\). Here, for each \(t\in \N\), \(x_t\) and \(u_t\) are the state and the control variable, and \(f(\cdot, \cdot)\) is a known mapping. Let \(U \Let \left(u_0, \ldots, u_{\horizon-1}\right)\), the MPC method entails solving a constrained optimal control problem (OCP) over a given time horizon \(\horizon \in \N\), formulated as
\begin{equation}\tag{P}
\label{eq:MPC:intro}
\begin{aligned}
& \inf_{U}  && \sum_{t=0}^{\horizon-1} \cost\bigl(\st_t, \ut_t\bigr) + \fcost\bigl(\st_{\horizon}\bigr) \nn \\
& \sbjto  && 
\begin{cases}
\st_{t+1}= \field(\st_t, \ut_t), \; \st_0 \Let \xz,\,\st_t \in \Mbb,\\
 \ut_t \in \Ubb \text{ for }t \in \aset[]{0,\ldots,\horizon-1},\, 
\st_{\horizon} \in \admfinst.
\end{cases}
\end{aligned}
\end{equation}
In this context, \(\cost(\cdot,\cdot)\) and \(\fcost(\cdot)\) denote the stage cost and the final (or terminal) cost functions, respectively, while \(\Mbb\), \(\Ubb\), and \(\admfinst\) define the admissible sets for states, control inputs, and final states, respectively. At \(t = 0\), the system's initial state is observed as \(\st_0 \Let \xz\), and the optimization problem is solved with \(\xz\) to determine an optimal control sequence \(U\as \in \admact^{\horizon}\); notice that \(\xz\) enters \eqref{eq:MPC:intro} as a \emph{parameter}. The first control input, \(u\as_0\), is then applied to the system. Time is advanced to \(t = 1\), and the process is repeated iteratively. This iterative approach inherently defines a \emph{feedback map} (potentially set-valued), as the control action \(\xz \mapsto u\as_0(\xz)\) implicitly depends on the current state \(\xz\). The goal of \embf{Explicit MPC} (ExMPC) is to computationally extract this \emph{implicit} map \cite{ref:weinan2022empowering}. The ExMPC approach is particularly valuable as it replaces the need for repeated numerical optimization in the MPC problem with a simple function evaluation at each state. This eliminates the necessity of solving the OCP \eqref{eq:MPC:intro} at every time step via shifting the time horizon, significantly reducing computational overhead for real-time applications. Consequently, ExMPC is well-suited for a variety of industrial process applications \cite{ref:Morari_explicit}, \cite{ref:MPC_for_power_elec}. Most existing ExMPC methods rely on multiparametric programming \cite[Chapter 7]{ref:XiLi-19}, but their applicability is largely confined to nominal linear, low-dimensional systems with short time horizons. For systems with moderate to high dimensionality, neural networks have been utilized to develop ExMPC algorithms, as seen in earlier works like \cite{ref:parsini1995} and more recent advancements such as \cite{ref:Mesbah:ExMPC,ref:weinan2022empowering,ref:ML:MultiParam:wang2024,ref:FastExplicitConvexNN,ref:karg:ExMPC,ref:ML:MPC:tutorial}. Although some of these approaches offer NN-based ExMPC methods with closed-loop guarantees, they all fail to adequately address the challenges of approximation errors or provide strategies for controlling them, even in low-dimensional settings.

While ExMPC, has become standard for nominal linear models, accounting for uncertainty in the systems introduces a fresh complication into the problem --- in particular, \eqref{eq:MPC:intro} becomes a minmax optimization problem. Robust optimization techniques can in principle be employed to ensure that the MPC policy becomes immune to a range of model uncertainties, but this is easier said than done due to the intrinsic difficulties associated with numerical solutions of robust optimization problems.
\subsection*{Our contributions}
\begin{enumerate}[label=\textup{\Alph*.}, leftmargin=*]
\item \label{contrib:aaMPC} \textbf{Approximation-aware formulation:} We present a numerically tractable technique to learn explicit feedback policies for a class of MPC problems with uniform approximation error guarantees. As a first step, one of our key contributions is the formulation of an \emph{approximation-aware} MPC setup that yields a policy that respects all constraints by design irrespective of approximation errors, thereby ensuring recursive feasibility and closed-loop stability. This is a clear point of departure from the existing literature that mostly provides providing a posteriori error bounds without taking such errors into account at the planning stage.
\begin{itemize}[leftmargin=*]
\item We build our theoretical results around a robust constrained MPC framework for an uncertain discrete-time linear dynamical system, where we explicitly account for the error arising due to the application of the approximate policy \(\xz \mapsto \mutrunc(\xz)\) instead of the (optimal) MPC policy \(\xz \mapsto \upopt(\xz)\), and robustify against that error. Here we introduce the \emph{approximation-aware} MPC problem: For this problem to be well-posed, it is essential to ensure a \emph{uniform} approximation error between the optimal policies \(\upopt(\cdot)\) and its approximants \(\mutrunc(\cdot)\). The algorithms established in the sequel inherently address this requirement.
        
\item The synthesis of the approximate policies relies on interpolation and learning techniques, and three attributes of this procedure are central:
\begin{enumerate}[label=\highlightg{(T-\alph*)}, align=left, widest=b, leftmargin=*]
\item \label{training1} the ability to guarantee uniform approximation within a preassigned error margin;
\item \label{training2} a training procedure that produces such an approximant, and
\item \label{training3} numerical viability.
\end{enumerate}        
We explore two learning engines in this context: quasi-interpolation and neural networks. For a prespecified uniform error tolerance \(\varepsilon > 0\), our quasi-interpolation driven procedure satisfies all the three attributes (\ref{training1}--\ref{training3}) for low-dimensional systems (at most \(\dimsp = 5\)); the neural network driven procedure checks the first attribute under abundance of training data but the second attribute is difficult to ensure in practice as we shall demonstrate via numerical experiments. There is, however, reason to stay cautiously optimistic concerning the applicability of neural networks, as we point out in the sequel.
\end{itemize}
\begin{figure}
    \centering
    \includegraphics[scale=0.35]{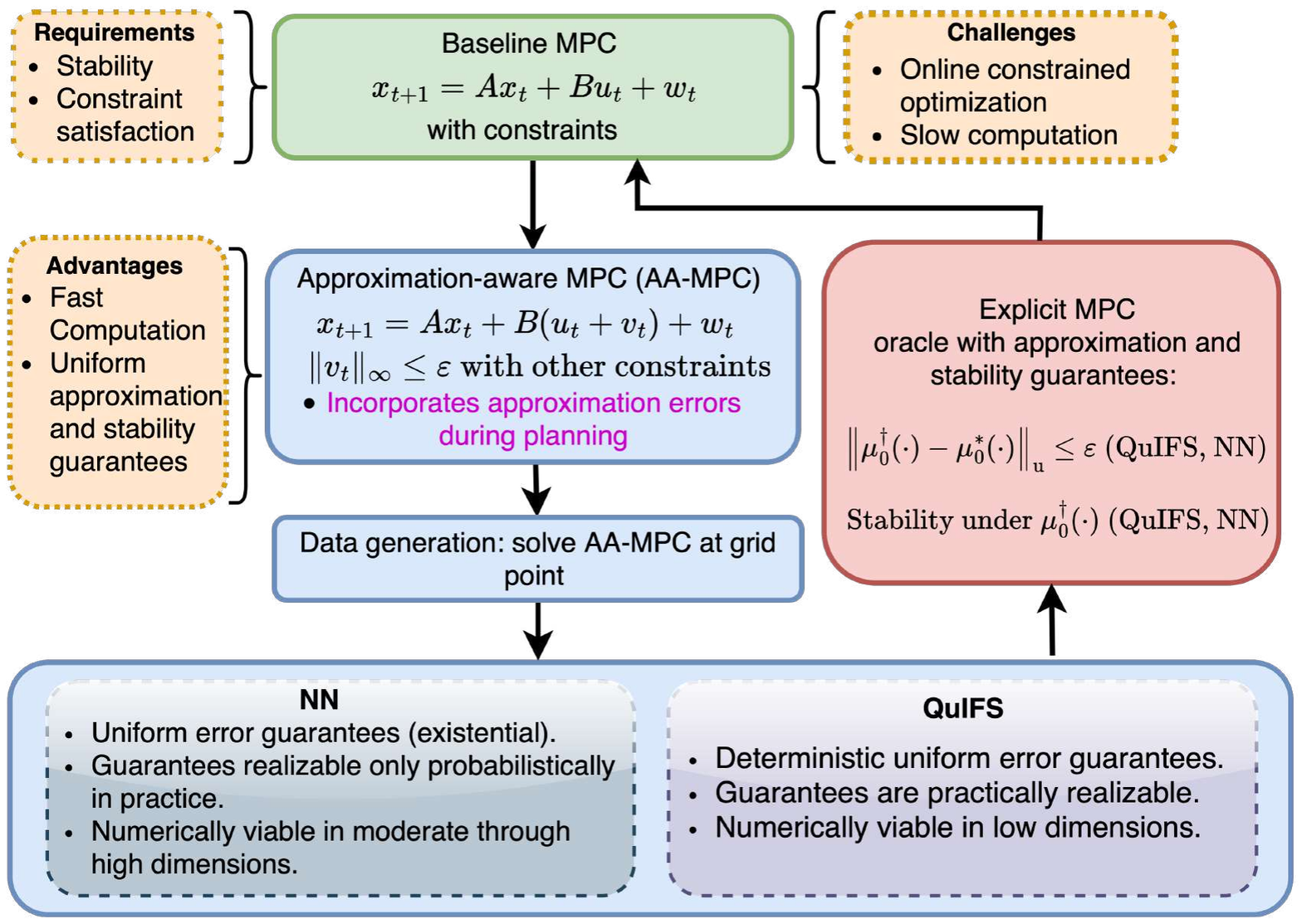}
    \caption{A bird's-eye view of our algorithm and it main components.}
    \label{fig:mcrf:schematic}
\end{figure}
\item \label{contrib:MSA} \textbf{A semi-infinite programming-based method for exact \texttt{(state, action)} data generation:} 
The underlying dynamical system in our framework is an \emph{uncertain linear system}, where the robust MPC problem is formulated as a worst-case minmax optimization. The (outer) minimization is performed over a class of control policies, while the (inner) maximization considers all admissible uncertainties, subject to convex constraints. The ensuing convex robust optimization problem is, in general, NP-hard. However, a recently-developed technique \cite{ref:DasAraCheCha-22} for \emph{exact} solutions to convex semi-infinite programs (which we refer to as the MSA algorithm in the sequel) are employed herein to solve this important class of robust MPC problems in an \emph{exact} fashion at specific grid points on the state space; we submit that this approach leads to the extraction of accurate optimizers for any feasible initial state.\footnote{See \S\ref{app:msa_overview} in Appendix \ref{sec:prelim} for a brief overview of the MSA algorithm.} Moreover, due to the absence of conservative approximations, \emph{the exact feasible set} of the minmax problem is larger than those obtained via conservative approximations; see Example \ref{exmp:rmpc_1} in \S \ref{sec:num_exp} for an illustration.\footnote{This example demonstrates that the feasible set of initial points \(\xz\) obtained using our method is strictly larger compared to the feasible set obtained by employing YALMIP's \cite{ref:YALMIP_lofberg2004} robust optimization module \cite{ref:lofberg2012automatic} which employs scenario-like approximation methods.} We draw attention to the fact that alternative techniques such as the scenario approach for solving the minmax problem may suffer, especially in high dimensions, from acute suboptimality issues highlighted in \cite[Section 1.5]{ref:MA:DC:RB} and \cite[Remark 12]{ref:DasAraCheCha-22}, which are circumvented via the MSA algorithm.\footnote{The MSA algorithm is not limited to quadratic costs and linear constraints; while certain specific structures of minmax problems may be solvable by alternate means, MSA offers a general framework that applies to all minmax problems of the ``convex'' type --- see \cite{ref:DasAraCheCha-22} for details.} Our results concerning this contribution are recorded in \S \ref{app:MSA}.

	The issue of computing the exact optimizers for the underlying minmax problem at each given initial state is addressed by the MSA algorithm, but the algorithm is typically slow (due to the presence of a global optimization step therein). Consequently, a direct online application may not be possible always, and that brings us to our second key contribution by means of function learning:

	\item \label{contrib:explicit} \textbf{Approximation/learning feedback maps from data:} Instead of computing control actions through optimization at each step, we learn a feedback map to evaluate them, but with a twist. In contrast to the standard approaches we advocate \emph{tight approximations} of the underlying implicit feedback for the construction of the explicit map; this approximation is driven by means of input-output data of the form \texttt{(state, action)}-pairs made available by the MSA algorithm. Only uniformly approximate designs can be utilized in this context for operational reasons (see Remark \ref{r:robust and uniform} for a detailed discussion), and we submit two approaches for learning such uniformly approximate explicit feedback maps:

\begin{enumerate}[leftmargin=*, widest=ii, align=left]
\item \label{contrib:explicit:qi} Along the lines of the \(\textsf{QuIFS}\) algorithm reported in \cite{ref:GanCha-22}, the first approach applies to low-dimensional systems and employs \emph{quasi-interpolation} to learn continuous functions from the \texttt{(state, action)} data on a uniform cardinal grid up to a \emph{preassigned uniform error}. For engineers, this paradigm is extremely well-suited, and works well in the robust MPC context. The designer specifies an error tolerance \(\varepsilon > 0\) and selects a rapidly decaying smooth generating/basis function \(\psi(\cdot)\) along with three parameters \((h, \mathcal{D}, \rzero) \in \loro{0}{+\infty}^3\) such that the error between the optimal and learned feedback policies remains within \(\varepsilon\) \emph{uniformly} over the admissible set. The \emph{uniformity} in the estimate permits us to accommodate, with foresight, this (small) error into the minmax formulation itself so that recursive feasibility is not compromised at any stage and the closed-loop robust stability is ensured. Our results corresponding to this approach are reported in \S\ref{subsec:main_first_approach} -- \S\ref{sec:main_tech_results}. 

\item \label{contrib:explicit:nn} The presence of a uniform cardinal grid on the feasible set in \ref{contrib:explicit:qi} limits this interpolation-based technique to low state dimensions. For moderate-dimensional systems, recent results regarding uniform approximation properties of deep neural networks \cite{ref:Shen_NN_2022} point to an alternative approach to learn an approximate (and sometimes theoretically the \emph{exact}) explicit feedback map based on the offline precomputed \texttt{(state, action)} data (in this regime either on a uniform cardinal grid or on a suitably randomly sampled set of points as the case may be) via training an appropriate neural network (in our case, a deep ReLU network \cite{devore_hanin_petrova_2021}). To fulfill the requirement \ref{training1}--\ref{training2} numerically, our empirical results demonstrate that there are reasons to be cautiously optimistic even though the \emph{current} state-of-the-art theoretical requirements are prohibitive from a numerical tractability standpoint. The results corresponding to this approach are reported in \S\ref{subsec:main_second_approach} -- \S\ref{sec:main_tech_results_nn}.
\end{enumerate}
\end{enumerate}
\textbf{Notation:} We let \(\N \Let \aset{1,2,\ldots}\) denote the set of positive integers, \(\Nz\Let \Nzr\),  and let \(\Z\) denote the integers. Let \(\hat{a},\ol{a} \in \N\) and \(\ol{a}>\hat{a}\), then \([\hat{a};\ol{a}]\Let \aset[]{\hat{a},\hat{a}+1,\ldots,\ol{a}}\). The vector space \(\Rbb^d\) is assumed to be equipped with standard inner product \(\inprod{v}{v'}\Let \sum_{j=1}^d v_j v'_j\) for every \(v,v' \in \Rbb^d\). For any arbitrary subset \(X \subset \Rbb^d\) we denote the interior of \(X\) by \(\intr X\). The notation \(\oplus\) is the standard Minkowski sum of two sets. We denote the uniform function norm via the notation \(\unifnorm{\cdot}\); more precisely, for a real-valued bounded function \(f(\cdot)\) defined on a set \(S\), it is given by \(\unifnorm{f(\cdot)} \Let \sup_{x \in S}|f(x)|\). For vectors residing in some finite dimensional vector space \(\Rbb^d\) we employ the notation \(\norm{\cdot}_{\infty}\) for uniform norm. For \(\ell \in \N\) and for \(p \in \lcrc{1}{+\infty}\), the \(\ell\)-dimensional closed ball centered at \(x\) and of radius \(\eta>0\) with respect to the \(p\)-norm (vector) is denoted by \(\Ball^{\ell}_{p}[x, \eta]\); \(\Gamma(\cdot)\) will denote the standard Gamma function. The Schwartz space of rapidly decaying \(\Rbb\)-valued functions \cite[Chapter 2]{ref:mazyabook} on \(\Rbb^d\) is denoted by \(\mathcal{S}(\Rbb^d)\).

\section{Problem formulation}\label{sec:prob_form}
Consider a time-invariant discrete-time control system
\begin{equation}
    \label{eq:system}
    \st_{t+1} = Ax_t+Bu_t+w_t,\quad \st_0=\xz\text{ given},\quad t \in \Nz,
\end{equation}
along with the following data:
\begin{enumerate}[label=\textup{(\ref{eq:system}-\alph*)}, leftmargin=*, widest=b, align=left]
\item \label{eq:system-st-con-dist}\(\st_{t} \in \Rbb^{d}\), \(\con_t \in \Rbb^m\), and \(\dist_t \in \Rbb^d\) are the vectors representing the states, control inputs and uncertainties at the time \(t\), with \(A \in \Rbb^{d \times d}\) and \(B \in \Rbb^{d \times m}\);

\item \label{eq:constraint-sets-st-ut} the constraints on the state and control inputs of the system \eqref{eq:system} are given by
\(
\st_t \in \Mbb \text{ and } \ut_t \in \Ubb \text{ for all }t \in \Nz,
\)
where \(\Mbb\subset \Rbb^d\) and \(\Ubb \subset \Rbb^m\) are nonempty compact and convex subsets both containing the origins \(0 \in \Rbb^d\) and \(0 \in \Rbb^m\) in their respective interiors. The vector \(\dist_t\) represents the uncertainty entering in the system at each time \(t\) and
\(\dist_t \in \Wbb\) for all \(t \in \Nz,\) where \(\Wbb \subset \Rbb^d\) is a nonempty compact and convex with origin \(0 \in \Rbb^d\) in its interior.
\end{enumerate}
We assume that the following objects are given:
\begin{enumerate}[label=\textup{(\ref{eq:system}-\alph*)}, leftmargin=*, widest=b, align=left, start=3]
    \item \label{eq:system:horizon} a time horizon \(\horizon \in \N\); a jointly continuous and convex \emph{cost-per-stage} function
		\(
		\cost(\cdot,\cdot)
		\)
		and a continuous and convex \emph{final-stage cost} function
		\(
			\fcost(\cdot)
		\)
		satisfying \(\cost(0,0)=0\) and \(\fcost(0)=0\);
  
    \item \label{eq:system:terminal_set} a closed and convex terminal set \(\admfinst \subset \Mbb\) containing the origin \(0 \in \Rbb^d\) in its interior; a class of admissible control policies \(\policies\) consisting of a sequence \( \policy(\cdot) \Let (\policy_t(\cdot))_{t=0}^{\horizon-1}\) of continuous maps such that \(\pi_i: \Mbb \lra \Ubb\) for each \(i\) and we write \(\dummyu_t = \pi_t(\dummyx_t)\).
\end{enumerate}
\vspace{1mm}
Let \(\winopt \Let (w_t)_{t=0}^{\horizon-1}\). Given the preceding set of data, the \textbf{baseline robust MPC problem (BL-MPC)} for the system \eqref{eq:system} accompanied by the data \ref{eq:system-st-con-dist}-\ref{eq:system:terminal_set} is:
\begin{equation}
	\label{eq:baseline robust MPC}
	\begin{aligned}
		& \inf_{\policy(\cdot)} \sup_{\winopt} && \mathbb{J}_{\horizon} \bigl(\xz,\pi(\cdot),\winopt \bigr) \Let \sum_{t=0}^{\horizon-1} \cost(\dummyx_t, \dummyu_t) + \fcost(\dummyx_{\horizon})\\
		& \sbjto && \hspace{-3mm}\begin{cases}
			\text{dynamics }\eqref{eq:system},\,\dummyx_0 = \xz,\,
			\dummyu_t = \policy_t(\dummyx_t),\\ \dummyx_\horizon\in\admfinst,\,
			\dummyx_t\in\Mbb \text{ and } \dummyu_t\in\admact\\ \text{for all }(t,w_t) \in [0;\horizon-1] \times \Wbb.
		\end{cases}
	\end{aligned}
\end{equation}
A solution to \eqref{eq:baseline robust MPC} (if exists) is a control policy \(\opt{\policy}(\cdot) \Let \bigl(\opt\policy_t(\cdot)\bigr)_{t=0}^{\horizon-1}\), which, by definition, satisfies all constraints for any disturbance realization. If \(\opt{\policy}(\cdot)\) is unique for each \(\xz \in \Mbb\) (a parameter of \eqref{eq:baseline robust MPC}), the mapping from \(\xz\) to \(\opt{\policy}(\cdot)\) is a function, with \(\opt{\policy}_0(\cdot)\) being the \emph{feedback map}. Non-uniqueness results in a set-valued mapping.\footnote{See \cite[Lecture 0, p. 4]{ref:dontchev_variational} for details.} Computing the exact \(\opt{\policy}_0(\cdot)\) is extremely challenging, as it requires solving the Bellman-Isaacs equation \cite{ref:TM:book:2008}. Instead, we adopt an optimization-based approach. Assuming \(\opt{\policy}(\cdot)\) can be evaluated at each feasible point, we construct a tight uniform approximation of \(\opt{\policy}_0(\cdot)\). Naturally, this introduces approximation errors, which manifest as uncertainties in the control actions and these must be accounted for at the design stage to ensure recursive feasibility, leading to the focus of the next subsection.

\subsection{From the (BL-MPC) \eqref{eq:baseline robust MPC} to the approximation-aware MPC problem}\label{subsec:BLMPC_to_ARMPC}
To account for approximation error, we solve an OCP that exhibits several common features with \eqref{eq:baseline robust MPC} at the level of the cost and the constraints, but differs at the level of the dynamics. To this end, define \(\agrdist_t \Let w_t+Bv_t\) which takes values in the set \(\agradmdist \Let \Wbb \oplus B\Vbb\) and let \(W \Let \bigl(\agrdist_0,\ldots,\agrdist_{N-1}\bigr) \in \agradmdist^{\horizon}\). Consider the dynamical system
\begin{align}\label{eq:noisy-system}
    \st_{t+1} = Ax_t+Bu_t+\agrdist_t,\quad x_0=\xz\text{ given},\quad t \in \Nz,
\end{align}
where \(x_t \in \Mbb\), \(u_t \in \Ubb\), and the uncertainly component \((w_t ,v_t) \in \Wbb \times \Vbb\) for each \(t=0,\ldots,\horizon-1\), where the designer-specified uniform approximation error margin \(\varepsilon > 0\) has been fixed \textit{a priori}. The data \ref{eq:system-st-con-dist}--\ref{eq:constraint-sets-st-ut} carry over to \eqref{eq:noisy-system} modulo obvious changes; and the data \ref{eq:system:horizon}--\ref{eq:system:terminal_set} from the baseline MPC problem \eqref{eq:baseline robust MPC} are satisfied. Our \textbf{approximation-aware robust MPC problem (AA-MPC)} for the synthesis of receding horizon feedback is given by
\begin{equation}
\begin{aligned}\label{eq:approx-ready robust MPC}
		& \inf_{\policy(\cdot)} \sup_{W} && \mathbb{J}_{\horizon} \bigl(\xz,\pi(\cdot),W \bigr) \Let \sum_{t=0}^{\horizon-1} \cost(\dummyx_t, \dummyu_t) + \fcost(\dummyx_{\horizon})  \\
		& \sbjto && \hspace{-3mm}\begin{cases}
			\text{dynamics }\eqref{eq:noisy-system},\,\dummyx_0 = \xz,\,
			\dummyu_t = \policy_t(\dummyx_t),\\\dummyx_\horizon\in\admfinst,\,
			\dummyx_t\in\Mbb \text{ and } \dummyu_t+ v_t\in\admact\\ \text{for all }(t,W) \in [0;\horizon-1] \times \agradmdist^{\horizon}.
		\end{cases}
\end{aligned}  
\end{equation}  
    
\begin{remark}\label{r:robust and uniform}
The noise \((v_t)_{t \in \Nz}\) in \eqref{eq:noisy-system} represents uncertainty in the control action introduced by approximating the optimal feedback in closed-loop. The constraint \(\|v_t\|_{\infty} \le \eps\) is enforced during the synthesis stage of \eqref{eq:approx-ready robust MPC}, ensuring all uncertainties in \eqref{eq:system} are accounted for. Our learning method depends on the state-space dimension. For low-dimensional systems, we use a grid-based quasi-interpolation technique \(\quifs\) (Quasi-interpolation-based Feedback Learning), while moderate-dimensional systems employ neural network-based approximation \(\nnfs\) (Neural Network-based Feedback Learning); both furnish uniform approximation guarantees, as detailed in \S\ref{sec:main_results}.
\end{remark}


\subsection{From the (AA-MPC) \eqref{eq:approx-ready robust MPC} to a convex SIP problem}\label{subsec:ARMPC_to_SIP}
The (BL-MPC) and the (AA-MPC), i.e., the robust OCPs \eqref{eq:baseline robust MPC} and \eqref{eq:approx-ready robust MPC}, in general, are numerically intractable despite being perfectly sensible in theory. For state-feedback parameterization of the control policy, it is well-known that the set of admissible decision variables can turn out to be \emph{non-convex} \cite{ref:Lof-03} in general. We adopt a parameterization employing the past disturbance data instead \cite{ref:Lof:minmax:cdc}
\begin{align}\label{eq:distfb_par}
    u_t = \sum_{i=0}^{t-1} \theta_{t,i}\agrdist_i +  \eta_t,\quad t=0,\ldots,\horizon-1,
\end{align}
where each \(\theta_{t,i} \in \Rbb^{m \times d}\), \(\eta_t \in \Rbb^m\), and eliminating \((t,i)\) we have \((\theta,\eta) \in \Rbb^{ m\horizon\times d\horizon } \times \Rbb^{m\horizon}\). Employing the  parameterization in \eqref{eq:distfb_par}, the problem \eqref{eq:approx-ready robust MPC} now reads as
\begin{equation}
\label{eq:approx ready param robust MPC}
\mathbb{J}\as_{\horizon}(\xz) \Let  \inf_{\theta,\eta} \sup_{W} \, \aset[\big]{\mathbb{J}_{\horizon}\bigl(\xz,\theta,\eta,W\bigr) 
\suchthat \text{the constraints in \eqref{eq:approx-ready robust MPC}}}.
\end{equation}
Let the set of feasible decision variables \((\theta,\eta) \in \Rbb^{ m\horizon\times d\horizon } \times \Rbb^{m\horizon}\) for the problem \eqref{eq:approx ready param robust MPC} is given by $\mathcal{A}_{\theta,\eta}(\xz)$. Notice that, for a given initial state \(\xz\), the set \(\mathcal{A}_{\theta,\eta}(\xz)\) is convex \cite{ref:Lof-03}. This approach allows for the determination of an admissible policy by solving a \emph{convex semi-infinite program}, for a wide range of disturbances and uncertainties.

We define the feasible set of initial states, corresponding to the optimization problem \eqref{eq:approx ready param robust MPC} by \begin{align}\label{feas_set_X_N}\fset\Let \aset[]{\xz \in \Rbb^d \mid \mathcal{A}_{\theta,\eta}(\xz) \neq \emptyset}.\end{align}

\begin{assumption}\label{assump:param_feasible_set}
We assume that the set of feasible initial states \(\fset\) defined in \eqref{feas_set_X_N} is nonempty and that the problem \eqref{eq:approx ready param robust MPC} admits a unique solution that we denote by \(\bigl(\mu_t\as(\cdot)\bigr)_{t=0}^{\horizon-1}.\)
\end{assumption}

The first entry \(\upopt(\cdot)\) of the optimal policy \(\bigl(\mu_t\as(\cdot)\bigr)_{t=0}^{\horizon-1}\) obtained by solving the OCP \eqref{eq:approx ready param robust MPC} is important for us, and for each \(\xz\), the evaluation \(\upopt(\xz)\) of \(\upopt(\cdot)\) at \(\xz\) is obtained numerically by solving \eqref{eq:approx ready param robust MPC}. As mentioned above, before the synthesis process begins, the designer is allowed to fix an error margin, say \(\varepsilon>0\), and, if \(\mutrunc(\cdot)\) is an approximation of \(\upopt(\cdot)\), then it is stipulated that the approximation error satisfies
\begin{equation}
	\label{eq:policy bound}
	\unifnorm{\upopt(\cdot) - \mutrunc(\cdot)} \le \varepsilon,
\end{equation}
i.e., \(\|\upopt(y)-\mutrunc(y)\| \le \varepsilon\) for all feasible \(y \in X_{\horizon}\).
Our \embf{learning-based robust feedback control} technique for \eqref{eq:system} is encoded in the following two steps:
\begin{enumerate}[label=\highlightg{(\textsf{FL}\arabic*)}, leftmargin=*, widest=b, align=left]
	\item \label{explicit_mpc_1} measure the states \(\st_t\) at time \(t\),
	\item \label{explicit_mpc_2} apply \(\con_t = \mutrunc(\st_t)\), increment \(t\) to \(t+1\), and repeat.
\end{enumerate}
\vspace{1mm}
Since the approximation error uncertainty \eqref{eq:policy bound} is accounted for in \eqref{eq:approx ready param robust MPC}, using \(\con_t = \mutrunc(\st_t)\) in \eqref{eq:system} ensures constraint satisfaction at every stage and recursive feasibility by design. Our learning-based robust feedback control will be implemented via steps \ref{explicit_mpc_1}--\ref{explicit_mpc_2} using the approximate policy \(\mutrunc(\cdot)\). For a streamlined presentation, we use the same notation, \(\mutrunc(\cdot)\), to represent the approximate policy learned through both quasi-interpolation and the neural network architecture.

We now transition towards the numerical technique that we devise herein to solve \eqref{eq:approx ready param robust MPC}. Consider the problem \eqref{eq:approx ready param robust MPC} with its problem data and notations. Introducing a \emph{slack variable} \(r \in [0,+\infty[\), the problem \eqref{eq:approx ready param robust MPC} can be translated into:
\begin{equation}
	\label{eq:approx ready param SIP robust MPC}
	\begin{aligned}
		& \inf_{\theta,\eta, r}  && r\\
		& \sbjto && \begin{cases}
		\mathbb{J}_{\horizon}\bigl(\xz,\theta,\eta,W\bigr)-r \le 0,\,
		r \in \lcro{0}{+\infty},\\
			\text{the constraints in \eqref{eq:approx-ready robust MPC}}, \text{ for all } W \in \agradmdist^{\horizon}
		\end{cases}
	\end{aligned}
\end{equation}
Observe that \eqref{eq:approx ready param SIP robust MPC} is a \embf{convex semi-infinite optimization problem} (CSIP) and it requires the constraints to be satisfied for an infinite number of admissible noise vectors \(W \in \agradmdist^{\horizon}\) and thus there are  (potentially) \emph{uncountably many} such constraints. We translated the problem \eqref{eq:approx ready param robust MPC} to the SIP \eqref{eq:approx ready param SIP robust MPC} in order to establish a tractable algorithm that solves \eqref{eq:approx ready param robust MPC} in an \embf{exact fashion} by means of the MSA algorithm (we refer to \S\ref{app:msa_overview} and \cite{ref:DasAraCheCha-22} for more details on the algorithm).

We direct our attention to the SIP \eqref{eq:approx ready param SIP robust MPC}. To this end, define 
\begin{align}\label{dvar}
\dvar \Let \dim(\theta)+\dim(\eta)+\dim(r), 
\end{align}
\(\totdvar \Let \dvar \cdot \horizon\), and the function \(\agradmdist^{\totdvar} \ni \mathcal{W}  \Let \bigl( \conin{1},\ldots,\conin{\dvar}\bigr)\mapsto \gfunc(\mathcal{W};\xz) \in \Rbb\) for each \(\xz \in \fset\) by
\begin{align}\label{eq: inner param SIP robust MPC}
     \mathcal{G}\bigl(\mathcal{W};\xz\bigr)  \Let \inf_{\theta, \eta, r}  \left\{r  \;\middle\vert\;  
    \begin{array}{@{}l@{}}
       \text{constraints in } \eqref{eq:approx ready param SIP robust MPC} \; \text{for all} \\
        \conin{i}\in\agradmdist^{\horizon},\,i \in [1;\dvar]
        \end{array}
        \right\}.
\end{align}
Observe that the problem \eqref{eq: inner param SIP robust MPC} is a \emph{relaxed} version of the SIP \eqref{eq:approx ready param SIP robust MPC}, where all the constraints in \eqref{eq: inner param SIP robust MPC} need to be satisfied for \emph{finite} number of indices \(\conin{i}\) for \(i \in [1;\dvar]\); and the quantity \(\dvar\) has a precise characterization given in \eqref{dvar}. This is in contrast to the case in \eqref{eq:approx ready param SIP robust MPC}, where all the constraints needed to be satisfied \emph{for all} \(W \in \agradmdist^{\horizon}\), i.e., for all \(\agrdist_t \in \agradmdist\) for each \(t \in [0;\horizon-1]\). We shall refer to the problem \eqref{eq: inner param SIP robust MPC} associated with the function \(\gfunc(\cdot;\xz)\), for each \(\xz\), as the \emph{inner problem} in the sequel.

Finally, to extract the exact solution to the SIP \eqref{eq:approx ready param SIP robust MPC} one needs to globally optimize \(\mathcal{W} \mapsto \gfunc(\mathcal{W},\xz)\) over \(\mathcal{W} \in \agradmdist^{\totdvar}\). Precisely, consider the maximization problem
\begin{equation}\label{eq: outer param SIP robust MPC 0}
\sup_{\mathcal{W}} \,\,\aset[\big]{\mathcal{G}\bigl(\mathcal{W};\xz\bigr)\suchthat  \mathcal{W}\in \agradmdist^{\totdvar}},
\end{equation}
which will be called the \emph{outer problem} thereof, and this construct furnishes a single numerical method to solve \eqref{eq:approx ready param SIP robust MPC} using the MSA algorithm described in \cite{ref:DasAraCheCha-22}. As emphasized earlier, our algorithm algorithm solves the SIP \eqref{eq:approx ready param SIP robust MPC} \emph{exactly} in the sense that the \emph{optimal value} of the original SIP \eqref{eq:approx ready param SIP robust MPC} and the global maximization problem \eqref{eq: outer param SIP robust MPC 0} are the same (see Theorem \ref{thrm:value_func_equality}).

\section{Main results: Theory and Algorithms}\label{sec:main_results}     
In this section, we present our main theoretical results, with proofs available in Appendix \ref{appen:msa:proofs} and Appendix \ref{appen:ExMPC:proofs}. First, we offer a high-level overview of our findings:
\begin{enumerate}[leftmargin=*, widest=b, align=left]
	\item \label{result_1} (\S \ref{app:MSA}) The primary object is the minmax robust MPC problem \eqref{eq:baseline robust MPC} which we, via several intermediate steps, converted into the SIP \eqref{eq:approx ready param SIP robust MPC}. Various system theoretic properties of \eqref{eq:approx ready param SIP robust MPC} under the application of the MSA algorithm are recorded in \S \ref{app:MSA} along with a global optimization algorithm to solve \eqref{eq: outer param SIP robust MPC 0}; see Algorithm \ref{alg:msap} (\(\exctsol\)). The primary result is Theorem \ref{thrm:value_func_equality}; proof is given in Appendix \ref{appen:msa:proofs}.

	\item \label{result_2} (\S\ref{subsec:main_first_approach}, \S\ref{sec:main_tech_results}, \S\ref{subsec:main_second_approach}, \S\ref{sec:main_tech_results_nn}) We then leverage Algorithm \ref{alg:msap} and establish two different feedback learning algorithms. For a given \(\varepsilon>0\), both the algorithms theoretically provide guarantees for uniform approximation of the unknown optimal policy \(\upopt(\cdot)\) via the approximate policy \(\mutrunc(\cdot)\), under certain regularity assumptions on \(\upopt(\cdot)\). 
    \begin{itemize}[leftmargin=*]
        \item (\(\quifs\)) Algorithm \ref{alg:extension_algo}: based on quasi-interpolation, for low-dimensional systems. The primary result is Theorem \ref{thrm:stability_main_result}; proof is given in Appendix \ref{appen:ExMPC:proofs}. 
        \item (\(\nnfs\)) Algorithm \ref{alg: nn_algo}: based on neural networks, for moderate-dimensional systems. The primary result is Theorem \ref{thrm:stability_main_result_nn}; proof is given in Appendix \ref{appen:ExMPC:proofs}. 
    \end{itemize}     
\end{enumerate}
\subsection{The application of the MSA to the minmax problem}\label{app:MSA}
We begin by recording one of our chief observations concerning the problems \eqref{eq:approx ready param robust MPC} and \eqref{eq:approx ready param SIP robust MPC}. Let us fix some notations. Let \(\sipfeasbset \Let \Rbb^{ m\horizon\times d\horizon } \times \Rbb^{m\horizon} \times \lcro{0}{+\infty}\); by \(x_{\theta,\eta}(t;\xz,W)\) we mean the solution to the recursion \eqref{eq:noisy-system} at time \(t\) with initial state \(\xz\) under the control \eqref{eq:distfb_par} and disturbance sequence \(W\). We will continue to denote the terminal state \(x_{\theta,\eta}(\horizon;\xz,W)\) by \(x_N\) in the sequel.
\begin{lemma}
    \label{lem:Aux_lem}
    Consider the OCP \eqref{eq:approx ready param robust MPC} and the SIP \eqref{eq:approx ready param SIP robust MPC} with their associated data and notations. Define the set of admissible \((\theta, \eta, r)\) by
    \begin{align*}
        \feasSip \Let \aset[\Big]{(\theta, \eta, r) \in \sipfeasbset\suchthat \text{the constraints of the SIP \eqref{eq:approx ready param SIP robust MPC} hold}}.
    \end{align*}
    Then \(\feasSip\) is closed and convex.
\end{lemma}

The next assumption is a technical requirement which states that \eqref{eq:approx ready param robust MPC} is strictly feasible for \(\conin{i} \in \agradmdist^{\horizon}\) for each \(i=1,\ldots,\dvar\). This in turn implies strict feasibility of the optimization problem \eqref{eq:approx ready param SIP robust MPC}, and consequently, that of the relaxed problem \eqref{eq: inner param SIP robust MPC}. 
\begin{assumption}
\label{assum:slater's}
We stipulate for the OCP \eqref{eq:approx ready param robust MPC} that the following Slater-type condition holds: there exists a nonempty open set \(O \subset   \intr \admst \times \intr \admcont\) such that 
\begin{align*}
        \bigcap_{(\conin{i})_{i=1}^{\dvar}} 
        \left\{\bigl(\theta,\eta,r\bigr)\;\middle\vert\; 
\begin{array}{@{}l@{}}
        \mathbb{J}_{\horizon}(\xz,\theta^i,\eta^i,\conin{i}) -r < 0,\,\st_0=\xz,\st_{\horizon} \in \intr \admfinst,\\
        \bigl(\st^i_{\theta,\eta}(t;\xz,\conin{i}),u_t^i+v_t^i\bigr) \in O,\, \conin{i} \in \agradmdist^{\horizon}
        \end{array}
        \right\}
    \end{align*}
is nonempty. 
\end{assumption}

\begin{theorem}\label{thrm:value_func_equality}
Consider the OCP \eqref{eq:approx ready param robust MPC} and the corresponding SIP \eqref{eq:approx ready param SIP robust MPC} together with their associated data and notations, and let the Assumption \ref{assum:slater's} hold. Fix \(\xz \in \fset\) and consider the global maximization problem \eqref{eq: outer param SIP robust MPC 0}. Then:
\begin{enumerate}[label=\textup{(\ref{thrm:value_func_equality}-\alph*)}, leftmargin=*, widest=b, align=left]
\item \label{thrm:value_func_equality_0}
\(\agradmdist^{\totdvar} \ni  \mathcal{W} \mapsto \gfunc(\mathcal{W}; \xz) \in \Rbb\) is upper semicontinuous for every \(\xz \in \fset\),
\item \label{thrm:value_func_equality_1} there exists an \(\dvar\)-tuple \(\mathcal{W}\as\) that solves \eqref{eq: outer param SIP robust MPC 0}, and 
\item \label{thrm:value_func_equality_2} \(\mathbb{J}\as_{\horizon}(\xz)=\mathcal{G}\bigl(\mathcal{W}\as;\xz\bigr)\) for all \(\xz \in \fset\).
\end{enumerate}
\end{theorem}
\begin{remark}\label{rem:on_reformulation}
Theorem \ref{thrm:value_func_equality} asserts that it is sufficient to consider only \(\dvar\)-many disturbance realizations \(\conin{i}\) for \(i \in \aset[]{1,\ldots,\dvar}\) for which the constraints of the problem \eqref{eq: inner param SIP robust MPC} need to be satisfied; then the resulting relaxed problem solves the original SIP \eqref{eq:approx ready param SIP robust MPC}. To obtain such an \(\dvar\mbox{-}\)tuple optimal point \(\mathcal{W}\as\), the maximization problem \eqref{eq: outer param SIP robust MPC 0} in Theorem \ref{thrm:value_func_equality} must be solved globally on \(\agradmdist^{\totdvar}\), and at the optimal point \(\mathcal{W}\as\) the value of the minmax problem \eqref{eq:approx ready param robust MPC} and \eqref{eq: outer param SIP robust MPC 0} are identical. 
\end{remark}

\textbf{A global optimization based algorithm to solve \eqref{eq:approx ready param robust MPC}:}
Algorithm \ref{alg:msap} systematically addresses the global maximization problem \eqref{eq: outer param SIP robust MPC 0} arising from numerically solving \eqref{eq:approx ready param SIP robust MPC}. While any suitable global optimizer can be used, we employ simulated annealing for its scalability with dimensionality but do not recommend it universally.
\begin{algorithm2e}[!ht]
    \DontPrintSemicolon
    \SetKwInOut{ini}{Initialize}
    \SetKwInOut{giv}{Data}
    \giv{Stopping criterion $\SC(\cdot)$, threshold for the stopping criterion \(\tau\), fix \(\xz \in \fset.\)}
    \ini{initialize the constraint indices \( \bigl(W_{\mathrm{init}}^1, W_{\mathrm{init}}^2, \ldots, W_{\mathrm{init}}^{\dvar}\bigr) \in \agradmdist^{\totdvar}\), initial guess for \(\mathcal{G}_{\text{max}}\), initial guess for the solution \(\overline{\theta}\) and \(\overline{\eta}\)}
   
    \While{$\SC(n) \leqslant \tau$}
        {
        Sample (via simulated-annealing) the disturbance constraint set \(\mathcal{W}^n \Let  \bigl(W^{n,1},W^{n,2}, \ldots, W^{n,{\dvar}}\bigr) \in\agradmdist^{\totdvar}\)
     
            Solve the \emph{inner-problem} and evaluate \(\mathcal{G}_n = \mathcal{G}\bigl(\mathcal{W}^n;\xz\bigr) \) as defined in \eqref{eq: inner param SIP robust MPC}
            

            \emph{Recover} the solution  \(\theta_n,\eta_n,r_n \in \text{arg min}_{(\widetilde{\theta},\widetilde{\eta},\widetilde{r})} \aset[\big]{\widetilde{r} \,|\,\text{ constraints in } \eqref{eq: inner param SIP robust MPC} \text{ hold at } \pnt^n} \)
            
        \If{\(\mathcal{G}_n \ge \mathcal{G}_{\mathrm{max}}\)}
        {
            Set \( \mathcal{G}_{\text{max}}\) \(\longleftarrow\) \(\mathcal{G}_n\), \(\overline{\theta}\) \(\longleftarrow \) \(\theta_n\), \(\overline{\eta}\) \(\longleftarrow\) \(\eta_n\)  
            
             Update $n \gets n+1$ \;
        }
        }
   \caption{(\(\exctsol\)) Global optimization algorithm to solve \eqref{eq: outer param SIP robust MPC 0}}
\label{alg:msap}
\end{algorithm2e}

\subsection{\(\quifs\): feedback learning for low-dimensional systems}\label{subsec:main_first_approach}
We present an algorithm to learn an approximate policy \(\mutrunc(\cdot)\), offline, for low-dimensional systems, which, by design, satisfies the approximation bound \(\bigl\|\upopt(\cdot)-\mutrunc(\cdot)\bigr\|_{\mathrm{u}} \le \varepsilon\), for any preassigned \(\varepsilon > 0\). Our approximation scheme relies on the following assumption regarding the regularity of \(\upopt(\cdot)\).
\begin{assumption}\label{assum:LipCon_of_mu}
    We assume that the policy \(\xz \mapsto \upopt(\xz)\) is Lipschitz continuous with rank \(L_0\). 
\end{assumption}
\textbf{Algorithmic steps to generate \(\mutrunc(\cdot)\) via \(\quifs\):} We employ Algorithm \ref{alg:msap} to generate the approximation-aware policy \(\upopt(\cdot)\) at grid points of a uniform grid on the state-space \(\Mbb\). Armed with this \(\texttt{(state, action)}\) data \(\bigl(mh,\upopt(mh)\bigr)\) where \(mh \in \Mbb \,\cap\, \aset[]{mh \suchthat m\in \Z^d, \,h>0}\), before proceeding to learn \(\mutrunc(\cdot)\), we perform an extension procedure on the domain of \(\upopt(\cdot)\).\footnote{The extension is essential. Indeed, by design, \(\upopt(\cdot)\) is defined on \(X_N\); but to prove Theorem \ref{thrm:stability_main_result} ahead the domain of \(\upopt(\cdot)\) must be extended to \(\Rbb^d\). This is a technical requirement.} Moreover, we require the extended policy to maintain Lipschitz continuity.

To this end, let \(\widehat{X}_N \Let \fset \cap \aset[]{mh \suchthat m \in \Z^d}\) be an \(h\)-net \cite{ref:hnets} of \(X_N\); define the (component-wise) extended policy by 
\begin{align}\label{eq:extented_policy}
   \Rbb^d \ni x \mapsto \mu_{E}\as(x) \Let \inf_{y \in \widehat{X}_N} \bigl( \upopt(y) + L_0 \norm{x-y}_2 \bigr).
\end{align}
It is easy to check that \(\mu\as_{E}(\cdot)\) is Lipschitz on \(\Rbb^d\) with the same rank \(L_0\); see \cite[Appendix 3, Lemma 1]{ref:Loja_real}. In the sequel, we shall overload notation and continue to label the policy \(\mu_E\as(\cdot)\) after extension to \(\Rbb^d\) as \(\upopt(\cdot)\) itself. 

We are now ready to apply our first approximation scheme, based on quasi-interpolation, to learn \(\mutrunc(\cdot)\).\footnote{For a brief background on the quasi-interpolation scheme we refer to \cite[\S III]{ref:GanCha-22}.} For the extended policy \(\upopt(\cdot)\) the quasi-interpolation scheme is given by 
\begin{align}\label{eq:approx_policy_parent_mainres}
        \apprfb(\st) \Let \Dd^{-d/2}   \sum_{\mathclap{m \in \Z^d}}\upopt(mh)\,\psi\left(\frac{\st-mh}{h \sqrt{\Dd}}\right)
\end{align}
for all \(x \in \Rbb^d\); where \(\psi(\cdot)\) is a continuous generating/basis function that needs to satisfy 
\begin{itemize}[ leftmargin=*, widest=b, align=left]
\item \emph{the continuous moment condition of order \(M \in \N\)}, i.e.,
\begin{align}\label{eq:moment_condition}
	&\int_{\Rbb^d}\psi(y)\,\dd y=1\text{ and } \int_{\Rbb^d}y^{\alpha}\psi(y)\,\dd y=0 \nn\\& \,\,\text{for all}\,\, \alpha,\,1\le [\alpha] < M;
\end{align}
where \(\alpha\in \Nz^d\) is a multi-index and \([\alpha]=\alpha_1+\cdots+\alpha_d\). 
\item \emph{the decay condition}: For all \(\alpha \in \mathbb{N}^d\) satisfying \(0 \le [\alpha] \le \lfloor d/2 \rfloor + 1\), the function \(\psi(\cdot) \) is said to satisfy the decay condition of exponent \(K\) if there exist \(C_0>0\) and \(K>d\) such that 
\begin{align}\label{eq:decay_condition}
   \left(1+\|x\|\right)^K  |\partial^{\alpha} \psi(x)| \le C_0 \quad\text{for}\,\, x\in \Rbb^d.
\end{align}
\end{itemize}
Then from \cite[Theorem 2.25]{ref:mazyabook} or \cite[Theorem III.1]{ref:GanCha-22} we have the estimate 
\begin{align}\label{eq:lipschitz_estimate_mpc_proof_mainres}
    	\unifnorm{\apprfb(\cdot) - \upopt(\cdot)} \le C_{\gamma}L_{0}h\sqrt{\Dd}+ \Delta_0(\psi,\Dd) \,\,\text{on}\,\,\Rbb^d,
    \end{align}
where \(\Delta_0(\psi,\Dd) \Let \mathcal{E}_0(\psi,\Dd)\unifnorm{\upopt(\cdot)}\) is the saturation error, \(C_{\gamma} \Let M \cdot \Gamma (M)/\Gamma(M+2)\) is a constant (where \(\Gamma\) denotes the standard Gamma function), and the term
\(\mathcal{E}_0(\psi,\Dd)\) is 
\begin{align}\label{eq:saturation_decay_term}
    \mathcal{E}_0(\psi,\Dd)(\cdot) \Let \sup_{x \in \Rbb^d}\sum_{\nu \in \Z^d\setminus \{0\}}\mathcal{F} \psi(\sqrt{\Dd}\nu)\epower{2\pi i \inprod{x}{\nu}}.
\end{align}
Note that, for each \(x \in \Rbb^d\), \eqref{eq:approx_policy_parent_mainres} constitute a summation over the lattice \(\Z^d\). Instead, we use the finite-sum truncated approximant to learn our approximate policy. Define \( \finset_x(\rzero) \Let \aset[\big]{mh \suchthat mh \in \Z^d} \cap \Ball_2^d[x,\rzero h]\) and consider the sum
\begin{align}\label{eq:approx_policy}
 \mutrunc(\st) \Let \Dd^{-d/2}  \sum_{mh \in \mathbb{F}_x(\rzero)}\upopt(mh)\,\psi \left(\frac{\st-mh}{h \sqrt{\Dd}}\right)
\end{align}
for all \(x \in \Rbb^d\). Note that for any given \(\ol{\varepsilon} > 0\), by appropriately choosing \(\rzero\), the difference between \(\mutrunc(\cdot)\) and \(\mu_0^{\as}(\cdot)\) can be uniformly bounded by \(\ol{\varepsilon}\) for all \(x \in \mathbb{R}^d\). Hence, the use of \(\mutrunc(\cdot)\) is justified; see \cite[\S 2.3.3]{ref:mazyabook} for more details. 

Thus, three quantities in \eqref{eq:approx_policy} --- \(\rzero\), \(h\), and \(\Dd\) --- can be picked  \textbf{depending on the prescribed error margin} to ensure that the total \embf{uniform} error, stays within the preassigned bound \(\varepsilon\). The formal selection process for the tuple \((h, \Dd, \rzero) \in \loro{0}{+\infty}^3\) is outlined in Theorem \ref{thrm:stability_main_result} and Algorithm \ref{alg:extension_algo}.
\subsection{Approximation and stability guarantees under \(\quifs\)} \label{sec:main_tech_results}

\begin{assumption}\label{assump:stability_assumptions}
The following standard assumptions \cite{ref:MayFal-19} must be satisfied for recursive feasibility of the problem \eqref{eq:approx ready param robust MPC} and closed-loop stability of \eqref{eq:noisy-system} under the policy \(\upopt(\cdot)\):
\begin{itemize}[label=\(\circ\), leftmargin=*]
\item \label{eq:stability_prop1} There exists a continuous map \(\admfinst \ni \dummyx \mapsto \mu_F (\dummyx)\in \Ubb\) such that \(\fcost \bigl( A\dummyx+B\mu_{F}(x)\bigr) - \fcost(\dummyx) \le -\cost \bigl(\dummyx,\mu_{F}(\dummyx)\bigr)\) and \(A\dummyx+B\mu_F(\dummyx) + \agrdist  \in \admfinst\) for every \((\dummyx,\agrdist) \in \admfinst \times \agradmdist\).

\item\label{eq:stability_prop2} There exists \(\rho_1(\cdot) \in \classK\) such that for every \(\dummyx \in \admfinst\) and for every \(\agrdist\in \agradmdist\), \(\fcost(\cdot)\) satisfies
\(\fcost \bigl( A\dummyx+B\mu_F(\dummyx)+ \agrdist \bigr) -\fcost(\dummyx) \le - \cost\bigl(\dummyx,\mu_{F}(\dummyx)\bigr)+\rho_1 \bigl(\|\agrdist\|_{\infty}\bigr).\)\footnote{Definition of class \(\mathcal{K}\), \(\mathcal{K}_{\infty}\), and \(\mathcal{KL}\) functions can be found in \cite[Chapter 4, Definition 4.2]{ref:khalil}.}
\end{itemize}
\end{assumption}
\begin{remark}\label{rem:stability_value_bounds}
Given Assumptions \ref{assump:param_feasible_set} and \ref{assump:stability_assumptions}, it follows from \cite[\S3]{ref:MayFal-19} that the problem \eqref{eq:approx ready param robust MPC} remains recursively feasible. Furthermore, when \(\upopt(\cdot)\) is applied, the value function \(\valuefunc(\cdot)\) exhibits a \emph{descent behavior}
\begin{align}
	\label{eq:value_descent}
	\valuefunc(x{\nextst}) - \valuefunc(x) \leq -\cost(x, \upopt(x)) + \rho_1(W_{\text{sup}}).
\end{align}
for all \(x \in \fset\) and \(\agrdist \in \agradmdist\), where the successor state is defined as \(x{\nextst} \Let Ax + B \upopt(x) + \agrdist\), and the term \(W_{\text{sup}}\) represents the supremum of the disturbance magnitudes, i.e., \(W_{\text{sup}} \Let \sup \{ |\agrdist| \mid \agrdist \in \agradmdist \}\).
\end{remark}
Against this backdrop, here is our second technical result concerning the approximate feedback policy \(\mutrunc(\cdot)\) learned by \(\quifs\) (proof given in Appendix \ref{appen:ExMPC:proofs}):

\begin{theorem}\label{thrm:stability_main_result}
Consider the OCP \eqref{eq:approx ready param robust MPC} along with its associated data. Let \(\upopt(\cdot)\) be the corresponding approximation-aware MPC policy and suppose that assumptions \ref{assump:param_feasible_set} and \ref{assum:LipCon_of_mu} hold. Then:
\begin{enumerate}[label=\textup{(\ref{thrm:stability_main_result}-\alph*)}, leftmargin=*, widest=b, align=left]

\item \label{thrm:main:estimate} For every given \(\varepsilon>0\), there exist a generating function \(\psi(\cdot) \in \mathcal{S}(\Rbb^d)\), a triple \((h,\Dd,\rzero) \in \loro{0}{+\infty}^3\), such that the approximate feedback map \(\fset \ni \st \mapsto \mutrunc(\st) \in \Ubb\) defined in \eqref{eq:approx_policy} is within a uniform error margin \(\varepsilon\) from \(\upopt(\cdot)\); to wit,
		\begin{equation}
		    	\|\upopt(x) - \mutrunc(x)\| \le \varepsilon \quad \text{for all }\,x \in \fset.  \nn
		\end{equation}
\item \label{thrm:main:stability} Suppose that Assumption \ref{assump:stability_assumptions} holds and given the initial state \(x\), let \(x^{+}\) be the state at the next time instant. Then under \(\fset \ni \st \mapsto \mutrunc(\st)\) learned via Algorithm \ref{alg:extension_algo}: \(\quifs\), the closed-loop process \(x^{+} = Ax+B \mutrunc(x)+ w\) corresponding to the system \eqref{eq:system} is ISS-like stable in the sense of Definition \ref{def:isps_def}.
    \end{enumerate}
    \end{theorem}
    
\begin{algorithm2e}[!h]
\DontPrintSemicolon
\SetKwInOut{ini}{Initialize}
\SetKwInOut{giv}{Data}
\SetKwInOut{ext}{Extend}
\SetKwInOut{interpol}{Learning}
\giv{\(\upopt(\cdot)\) on \(\widehat{X}_N\)}
\ini{Lipschitz constant \(L_{0}\) of the policy \(\upopt(\cdot)\)}

\ext{Extend \(\upopt(\cdot)\) to \(\Rbb^d\) using \eqref{eq:extented_policy}}

\interpol{\(\circ\) Fix \(\varepsilon>0\) and choose the tuple \(\bigl(\psi(\cdot),h,\Dd,\rzero\bigr) \in \mathcal{S}(\Rbb^d) \times \loro{0}{+\infty}^3\); \\
		\(\circ\) Compute \(\mutrunc(\cdot)\) via \eqref{eq:approx_policy} and restrict \(\mutrunc(\cdot)\) to \(\fset\).} 
\caption{\(\quifs\): Quasi-interpolation-driven feedback synthesis}
\label{alg:extension_algo}
\end{algorithm2e}
\begin{remark}
Algorithm \ref{alg:extension_algo} represents a significant advancement over \cite{ref:GanCha-22}, as it satisfies all the properties \ref{training1}–\ref{training3} outlined in contribution (see \S\ref{s:intro}). It introduces a novel semi-infinite optimization-driven approach --- both theoretical and algorithmic --- for solving the minmax MPC problem \emph{exactly} while integrating it with a quasi-interpolation-based oracle for low-dimensional systems. In contrast, \cite{ref:GanCha-22} relied on off-the-shelf solvers for minmax optimization problems without contributing new methodological developments in this area.
\end{remark}
\subsection{\(\nnfs\): feedback learning for moderate-dimensional systems} \label{subsec:main_second_approach} 
Due to the use of a uniform grid, the synthesis approach established in \S\ref{subsec:main_first_approach}--\S\ref{sec:main_tech_results} is suitable for learning approximate policies only if the underlying dynamical system is of low dimension. To address this, particularly for moderate-dimensional systems, we develop an alternative approach using deep neural networks with ReLU activation functions. By leveraging recent results on uniform approximation using ReLU networks \cite{ref:Shen_NN_2022}, we provide guarantees for \emph{uniform} approximation error associated with the learning process, along with closed-loop stability guarantees, similar to Theorem \ref{thrm:stability_main_result}. See Theorem \ref{thrm:stability_main_result_nn} ahead in \S\ref{sec:main_tech_results_nn}.
    
\textbf{Algorithmic steps to generate \(\mutrunc(\cdot)\) via \(\nnfs\):}
Recall that \(\upopt(\cdot)\) is the approximation-aware policy generated via Algorithm \ref{alg:msap}. For learning the a feedback we will employ a feedforward ReLU neural network architecture \cite{devore_hanin_petrova_2021} and we denote the approximate policy to \(\upopt(\cdot)\) generated via the ReLU network by
\begin{equation}
    \label{eq:nn_approx_policy}
    \mutrunc(\cdot) \Let \Upsilon^{W,L}_{\upopt(\cdot)}(\relu(\cdot);d,d'),
\end{equation}
where for fixed \(W\) (width), \(L\) (depth), \(d\) (input dimension), \(d'\) (output dimension), and with \(x\mapsto \relu(x)\Let \max\,\{0,x\}\) (activation function), the function \(\Upsilon^{W,L}_{\upopt(\cdot)}(\relu(\cdot);d,d')\) denotes the ReLU approximation of the policy \(\upopt(\cdot)\). The \(\texttt{(state, action)}\) data for training the ReLU network is generated by sampling initial states \(\xz\) from the set \(\Mbb\) and solving \eqref{eq:approx ready param robust MPC} using the global optimization problem \eqref{eq: outer param SIP robust MPC 0} via Algorithm \ref{alg:msap}. We denote by \(\mathcal{X}_0\) the set of sampled initial states, obtained either through structured \emph{grid-based} sampling or \emph{uniform random} sampling over the domain \(\Mbb\). While a grid-based sampling is necessary for the quasi-interpolation-based construction of the approximate policy \(\mutrunc(\cdot)\) described in \S\ref{subsec:main_first_approach}, neural network approximation does not impose such restrictions, making the algorithm potentially scalable to systems with moderately many dimensions.

Recall Corollary \ref{corr:NN_approx_result_II} in Appendix \ref{appen:NN_results}. In our context, we replace \(f(\cdot)\) and \(f_{\relu}(\cdot)\) in Corollary \ref{corr:NN_approx_result_II} with \(\upopt(\cdot)\) and its ReLU-approximant \(\mutrunc(\cdot)\), respectively, and the estimate \eqref{eq:NN_estimate_2_Lip} can be rewritten as:
\begin{align}\label{eq:nn_approx_error_mainres}
     \unifnorm{\upopt(\cdot) - \mutrunc(\cdot)} &\le 131\sqrt{d} L_0 \bigg{(}\bigl(\tfrac{W}{3^{d+5}d}\bigr)^2 \bigl( \tfrac{L-18-2d}{22}\bigr)^2  \log_3 \bigl( \tfrac{W}{3^{d+5}d}+2 \bigr) \bigg{)}^{-1/d},
    \end{align}
where \(L_0\) is the Lipschitz rank of the feedback policy \(\upopt(\cdot)\) and \(W,L\) are the \emph{width} and \emph{depth} of the network respectively. The values of \(W\) and \(L\) are to be selected \textbf{according to the prescribed error margin}. To this end, we fix an error margin \(\varepsilon >0\),  and proceed with the aforementioned selection as follows:
\begin{itemize}[leftmargin=*, label=\(\triangleright\)]
    \item \label{pt: choose_w_l} We may begin by either setting a value to the network's width \(W\) or depth \(L\). The reasons influencing this choice include the facts that a greater width increases the number of trainable parameters, thus requiring a larger number of data points, but a greater depth increases training complexity owing to higher-order nonlinearities.

    \item In accordance with the choice made in the previous step, if the value of the width \(W\) is set first, then the value for the depth \(L\)  can be calculated using \eqref{eq:nn_approx_error_mainres}. Similarly, if the depth \(L\) is set first, then the width \(W\) found using the estimate \eqref{eq:nn_approx_error_mainres}. These estimates will be discussed in detail in \S\ref{sec:main_tech_results_nn}.
\end{itemize}
The above procedure ensures that the error between \(\upopt(\cdot)\) and \(\mutrunc(\cdot)\) stays uniformly within the prescribed bound \(\varepsilon\). The complete process is given in Algorithm \ref{alg: nn_algo}.
\begin{algorithm2e}[!ht]
\DontPrintSemicolon
\SetKwInOut{ini}{Initialize}
\SetKwInOut{giv}{Data}
\SetKwInOut{interpol}{Learning}
\giv{\(\upopt(\cdot)\) at points obtained from Algorithm \ref{alg:msap}}
\ini{Lipschitz rank \(L_{0}\) of the policy \(\upopt(\cdot)\)}

\interpol{
\hspace{-5mm}\begin{itemize}[label=\(\circ\),leftmargin=*]
    \item  Fix \(\varepsilon>0\) followed by setting a value for either \emph{depth} \(L\) or the \emph{depth} \(W\) for the ReLU network; 
    \item For a fixed width, using \ref{thrm:main:estimate_L}, or for a fixed depth, using \ref{thrm:main:estimate_W}, estimate the depth or the width depending on the fixed \(\varepsilon\);
    \item Generate the training dataset using Algorithm \ref{alg:msap} on a set of initial states \(\mathcal{X}_0\), sampled from \(\Mbb\). 
    \item Tune hyperparameters and train the NN to generate \(\mutrunc(\cdot)\) such that \(\bigl\|\upopt(x)-\mutrunc(x)\bigr\| \le \varepsilon\) for all \(x \in \Mbb.\) 
\end{itemize}
\vspace{1mm}
}
            
\caption{\(\nnfs\):Neural Network-driven feedback synthesis}
\label{alg: nn_algo}
\end{algorithm2e}
\subsection{Approximation and stability guarantees under \(\nnfs\)} \label{sec:main_tech_results_nn}
Here is our key result on ReLU neural network approximation (proof given in Appendix \ref{appen:ExMPC:proofs}):
\begin{theorem}\label{thrm:stability_main_result_nn}
  Consider the constrained OCP \eqref{eq:approx ready param robust MPC} along with its associated data. Let \(\upopt(\cdot)\) be the corresponding approximation-aware MPC policy and suppose that assumptions \ref{assump:param_feasible_set} and \ref{assum:LipCon_of_mu} hold. Then:
\begin{enumerate}[label=\textup{(\ref{thrm:stability_main_result_nn}-\alph*)}, leftmargin=*, widest=b, align=left]
\item \label{thrm:main:estimate_L} For every given \(\varepsilon>0\) and a fixed width \(\mathcal{W}_0\geq 3^{d+4}d\), there exists \(L \in \N\), and an approximate feedback map \(\fset \ni \st \mapsto \mutrunc(\st) \in \Ubb\) realized by a ReLU neural network of width $\mathcal{W}_0$ and depth $L$ such that the uniform error estimate 
\begin{equation}
\bigl\|\upopt(x) - \mutrunc(x)\bigr\| \le \varepsilon \quad \text{holds for all }x \in \fset.  \nn
\end{equation}
      
\item \label{thrm:main:estimate_W} For every given \(\varepsilon>0\) and a fixed depth \(\mathcal{L}_0\geq 29+2d\), there exists \(W \in \N\), and an approximate feedback map \(\fset \ni \st \mapsto \mutrunc(\st) \in \Ubb\) realized by a ReLU neural network of depth $\mathcal{L}_0$ and width $W$ such that the uniform error estimate
\begin{equation}
\bigl\|\upopt(x) - \mutrunc(x)\bigr\| \le \varepsilon \quad \text{holds for all }x \in \fset.  \nn
\end{equation}
      
\item \label{thrm:main:stability_nn} Suppose that Assumption \ref{assump:stability_assumptions} holds and given the initial state \(x\) let \(x^{+}\) be the state at the next time instant. Then under \(\fset \ni \st \mapsto \mutrunc(\st)\) constructed via Algorithm \ref{alg: nn_algo}: \(\nnfs\), the closed-loop process \(x^{+} = Ax+B \mutrunc(x)+w\) corresponding to the system \eqref{eq:system} is ISS-like stable in the sense of Definition \ref{def:isps_def}.
\end{enumerate}
\end{theorem}

\begin{remark}
In practice, the width and depth of ReLU networks required by Theorem \ref{thrm:stability_main_result_nn} are often significantly larger than what is empirically necessary. This discrepancy stems from the fact that the bounds in \cite{ref:Shen_NN_2022} were derived for Lipschitz continuous functions, which can result in overly conservative estimates. Our numerical experiments demonstrate that much shallower networks are often sufficient to achieve accurate approximations, indicating that the theoretical requirements may not always align with practical realities. Moreover, the non-exactness inherent in typical training procedures underscores the importance of a post-hoc validation of the approximated feedback map to ensure its reliability and effectiveness, thereby providing practical guarantees for the approximation.  
\end{remark}

\section{Numerical experiments}\label{sec:num_exp}
In this section, we focus on two examples to illustrate the numerical fidelity of the three main algorithms established in the preceding sections. The first example is the more comprehensive one --- it considers a minmax OCP for a second-order system, allowing clear visualization of the feedback maps learned by Algorithms \ref{alg:extension_algo} and \ref{alg: nn_algo}. Additionally, by demonstrating a larger feasibility region, this example highlights why Algorithm \ref{alg:msap} outperforms standard scenario-driven robust optimization methods. For the fourth-order example, we use only Algorithm \ref{alg: nn_algo} to construct the state and control trajectories. For both problems, we employed Algorithm \ref{alg:msap} via solving the \emph{inner problem} using MOSEK \cite{ref:mosek} in conjunction with YALMIP \cite{ref:YALMIP_lofberg2004} in MATLAB and for the \emph{outer problem}, which is a global maximization, we employed simulated annealing. For the simulated annealing, we used MATLAB's global optimization toolbox; we kept the initial temperature \(10\), and the number of iterations \(1500\) with the exponential temperature decay of the form \(T_t = 0.99T_{t-1}\), where \(T_t\) is the annealing temperature at the \(t^{\mathrm{th}}\) iteration. We performed our numerical computations related to data generation on MATLAB 2019b using the parallel computation toolbox in an \(36\) core server with Intel(R) Xeon(R) CPU E\(5-2699\) v\(3\), \(4.30\) GHz with \(128\) Gigabyte of RAM. 
\begin{example}\label{exmp:rmpc_1}
Consider the dynamical system \cite{ref:Gao_minmax_TAC,ref:GanchaLCSS}:
\begin{equation}\label{eq:exmp_0_dyn}
\st_{t+1}= \begin{pmatrix}
 0.732 & -0.086 \\ 0.172 & 0.990 
\end{pmatrix} \st_t + \begin{pmatrix}
0.060 \\ 0.006
\end{pmatrix} \ut_t + \begin{pmatrix}
     0.3 & 0.4 \\ 0.2 & 0.15
\end{pmatrix} w_t.
\end{equation}
Fix \(\horizon=5\), and define the discrete-time OCP
    \begin{equation}
	\label{eq:RMPC_num_ex_0}
	\begin{aligned}
		& \inf_{\policy(\cdot)} \sup_{\winopt} && \inprod{x_N}{Px_N} + \sum_{t=0}^{\horizon-1} \inprod{\st_t}{Q\st_t}+\inprod{u_t}{Ru_t}  \\
		& \sbjto && \begin{cases}
 \eqref{eq:exmp_0_dyn},\,\dummyx_0=\xz,\,
          \dummyx_t \in \Mbb,\text{ and }\dummyu_t \in \Ubb  \\
  \text{for all }(t,\winopt)\in [0;\horizon-1]\times \Wbb^{\horizon},
		\end{cases}
	\end{aligned}
\end{equation}
where \(\Mbb \Let \aset[]{x = (\dummyx_1,\dummyx_2) \in \Rbb^2 \suchthat \|\dummyx\|_{\infty}\le 1.5}\), \(\Ubb \Let \aset[]{\dummyu \in \Rbb \suchthat |\dummyu|\le 2}\), and \(\Wbb \Let \aset[]{\dummyw \in \Rbb \suchthat \|\dummyw\|_{\infty} \le 0.05}\). The state weighting matrix \(Q \Let I_{2\times 2}\) the \(2\times 2\)-identity matrix, the control weighting matrix \(R \Let 0.01\), and \(P \Let \begin{pmatrix}
5.5461 & 4.9873 \\ 4.9873 & 10.4940
\end{pmatrix}\). We fixed \(\eps=0.03\),
which is a uniform approximation error margin. According to our formulation \eqref{eq:approx ready param robust MPC} we defined the set \(\agradmdist \Let \lcrc{-0.05}{0.05} \oplus B \lcrc{-0.03}{0.03}\) and \(\agrdist_t = w_t+ B v_t\) takes values in \(\agradmdist\). Recall that \(W \Let (\agrdist_0,\ldots,\agrdist_{\horizon-1})\); next we employed the parameterization \eqref{eq:distfb_par}, and arrived at the approximation-aware OCP
\begin{equation}
	\label{eq:approx_ready_RMPC_num_ex_0}
	\begin{aligned}
		&  \min_{(\theta,\eta)} \max_{{W}} && \inprod{x_N}{Px_N} + \sum_{t=0}^{\horizon-1} \inprod{\st_t}{Q\st_t} + \inprod{\cont_t}{R\cont_t} \\
		& \sbjto && \begin{cases}
			x_{t+1} = Ax_t+B u_t+Bv_t+w_t,\\ \st_0=\xz,\,
			\dummyx_t \in \Mbb,\text{ and }\dummyu_t +\dummyv_t \in \Ubb\\
			\text{for all } (t,W) \in [0;\horizon-1] \times \agradmdist^{\horizon}.
		\end{cases}
	\end{aligned}
\end{equation}
\textbf{Data generation via Algorithm \ref{alg:msap}:} We translated \eqref{eq:approx_ready_RMPC_num_ex_0} as a CSIP of the form \ref{eq:approx ready param SIP robust MPC}, where the constraints are needed to be satisfied for an infinite number disturbance realizations \(\agrdist_t\) for each \(t=0,\ldots,\horizon-1\). But Algorithm \ref{alg:msap} indicates that we need only \(\dvar\)-many constraints to be encoded in the algorithm, and for the problem \eqref{eq:approx_ready_RMPC_num_ex_0}, \(\dvar = \mathrm{dim}(\theta)+\mathrm{dim}(\eta)+\mathrm{dim}(r)=50+5+1=56\). 

Next, the state-space \(\lcrc{-1.5}{1.5}\times \lcrc{-1.5}{1.5}\) was equipped with a uniform cardinal grid corresponding to \(h=0.02\) leading to \(22801\) points and at each of these points we solved the SIPs independently (offline) via Algorithm \ref{alg:msap} to generate the \(\texttt{(state, action)}\) data. The preceding choice of \(h\) is \emph{not} arbitrary and will become evident shortly, in particular this choice of \(h\) depends on the error-margin \(\eps\). The convergence of the simulated annealing associated with the outer problem can be verified from the rightmost subfigure of Figure \ref{fig:roa_comparison}, where we see that the algorithm almost converges after \(10^3\) iterations, and the variations between multiple runs of the simulated annealing are insignificant. 

Before moving to the learning procedure via Algorithm \eqref{alg:extension_algo} and Algorithm \eqref{alg: nn_algo}, we compared the region of attraction (i.e., the feasible region \(X_N\) for \eqref{eq:approx_ready_RMPC_num_ex_0}) of the feedback \(\upopt(\cdot)\) obtained by using Algorithm \ref{alg:msap} with the robust optimization module in YALMIP, creating a uniform grid of step size \(h=0.02\). Observe that:
\begin{enumerate}[leftmargin=*]
    \item The feasible region \(\fset\) obtained by using Algorithm \ref{alg:msap} is substantially larger than the one obtained from YALMIP's robust optimization module \cite{ref:lofberg2012automatic}; see the leftmost subfigure of Figure \ref{fig:roa_comparison}. This is a desirable feature of any MPC algorithm; a bigger region of attraction implies the feasibility of the problem for a larger number of points. In particular, solving the problem \eqref{eq:approx_ready_RMPC_num_ex_0} on a uniform grid with step size \(h=0.02\), Algorithm 1 produced \(700\) infeasible points, while the method in \cite{ref:lofberg2012automatic} yielded \(1436\) points. This represents a \(51.25\%\) reduction in infeasible points with Algorithm \ref{alg:msap}.
    \item The state trajectories obtained by employing Algorithm \ref{alg:msap} and YALMIP's robust optimization module with the region of attraction corresponding to Algorithm \ref{alg:msap} in the background, follow almost the same path and converge to a set containing the origin \(0 \in \fset \subset \lcrc{-1.5}{1.5}\times \lcrc{-1.5}{1.5}\); see right-hand subfigure in Figure \ref{fig:roa_comparison}.

    \item The optimal value obtained from Algorithm \ref{alg:msap} consistently exceeded that of the algorithm in \cite{ref:YALMIP_lofberg2004}; see Table \ref{tab:opt:val:comp}. This is because Algorithm \ref{alg:msap} provides exact solutions, closely tracking the original optimal value.
\end{enumerate}
\begin{table}[ht]
\centering
\caption{Comparison of optimal values obtained by Algorithm \ref{alg:msap} and the algorithm in \cite{ref:lofberg2012automatic}. }
\vspace{1mm}
\begin{tblr}{l c c}
\hline[2pt]
\SetRow{azure9}
\(\xz\) &
Value (Alg. \ref{alg:msap}) &
     Value (Alg. \cite{ref:lofberg2012automatic})\\ 
\hline
\((0.5,0.5)^{\top}\) & 2.198 &  2.14  \\
    \((-1, 0.4)^{\top}\)  & 1.794 &  1.489  \\
    \((1,-0.9)^{\top}\) & 2.623 & 2.25  \\
     \((-0.3,-0.7)^{\top}\)  & 2.48 &  2.412  \\
\hline[2pt]
\end{tblr}
\label{tab:opt:val:comp}
\end{table}
\begin{figure}[h]
  \begin{subfigure}[b]{0.49\linewidth}
    \includegraphics[width=6cm,height=5cm]{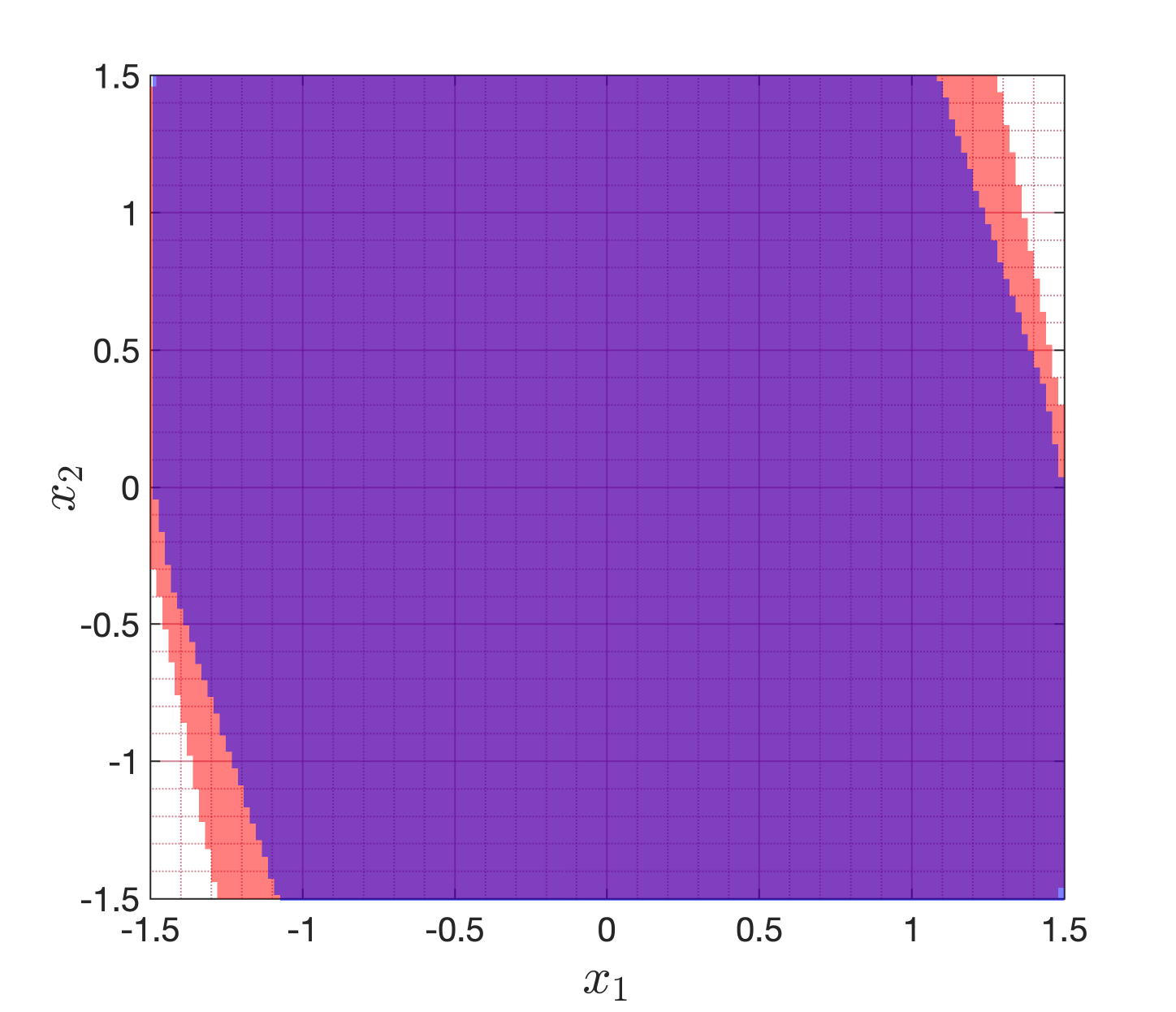}
  \end{subfigure}
  \begin{subfigure}[b]{0.49\linewidth}
    \includegraphics[width=6cm,height=5cm]{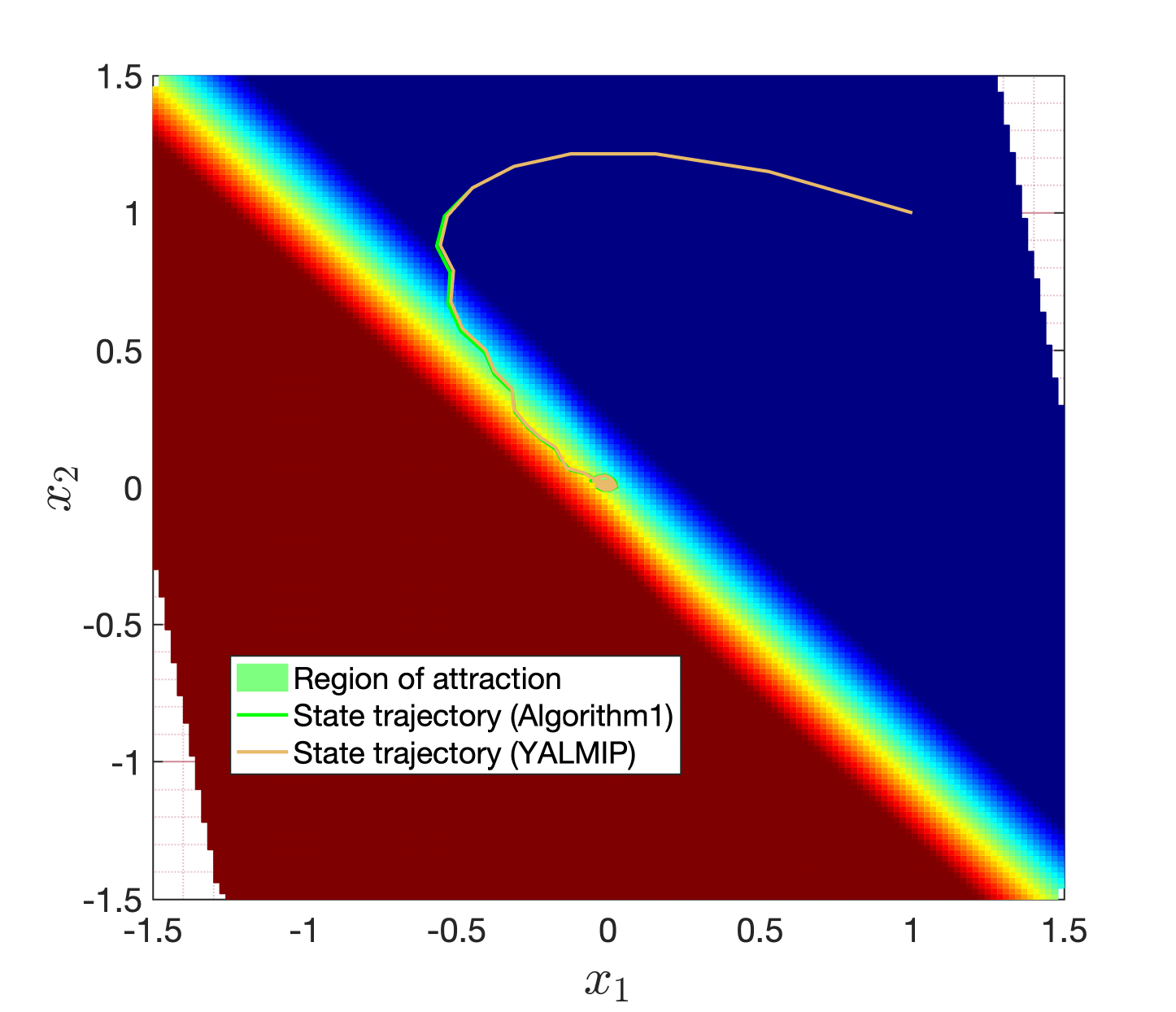}
  \end{subfigure}
\caption{The superimposed regions of attraction --- computed using YALMIP (violet) and Algorithm \ref{alg:msap} (magenta) --- are shown in the left-hand subfigure. The right-hand subfigure displays sample phase-space trajectories starting from the initial state \(\xz \Let (1\,\,1)^{\top}\), obtained using both Algorithm \ref{alg:msap} and YALMIP’s robust optimization module.}\label{fig:roa_comparison}
\end{figure}
\textbf{Learning via Algorithm \ref{alg:extension_algo}:} For the quasi-interpolation scheme, we employ the Laguerre-Gaussian basis function
\begin{align}\label{num:sixth_order_basis}
        x \mapsto \psi(x) \Let \frac{1}{\pi} \left( 3- 3 \|x\|^2+ \frac{1}{2}\|x\|^4\right)\epower{-\|x\|^2},
    \end{align}
that belongs to \(\mathcal{S}(\Rbb^2)\) and satisfies the moment condition of order \(M=6\), which means that the constant \(C_{\gamma} = 1/7\). Recall that we fixed the error margin \(\eps = 0.03\), i.e., we want \(\|v_t\|_{\infty} \le 0.03\). Next, fixing the shape parameter \(\Dd = 2\), we obtain the value of \(h = \frac{\varepsilon}{3C_\gamma L_0\sqrt{\Dd}} \approx 0.02\), where we have employed a rough, and conservative estimate of the Lipschitz constant \(L_0=2.5\) of \(\upopt(\cdot)\). We picked \(\rzero=4\). With these ingredients, after performing the extension procedure \eqref{eq:extented_policy}, for all \(x \in \Rbb^2\), we employed the quasi-interpolation formula
\begin{align}\label{num:sixth_quasi_interpolant}
     \mutrunc(x) \Let \frac{1}{\pi\Dd}\sum_{mh \in \finset_x(\rzero)}\upopt(mh)\biggl( 3 -  \frac{3 \|x-mh\|^2}{h^2\Dd}+\frac{2 \|x-mh\|^4}{h^4 \Dd^2}\biggr) \epower{-\frac{ \|x-mh\|^2}{\Dd h^2}}
 \end{align}
for the approximate feedback learning. The policies \(\mutrunc(\cdot)\) and \(\upopt(\cdot)\) can be seen in Figure \ref{fig:uopt_uapp_rmpc_ex_1} (the left-hand most and the right-hand subfigures, respectively), and the error surface is depicted in Figure \ref{fig:SipConvQError} (the right-hand subfigure), from where it is evident that the error tolerance margin is strictly satisfied. 
\begin{figure}[h]
  \begin{subfigure}[b]{0.49\linewidth}
    \includegraphics[width=6cm,height=5cm]{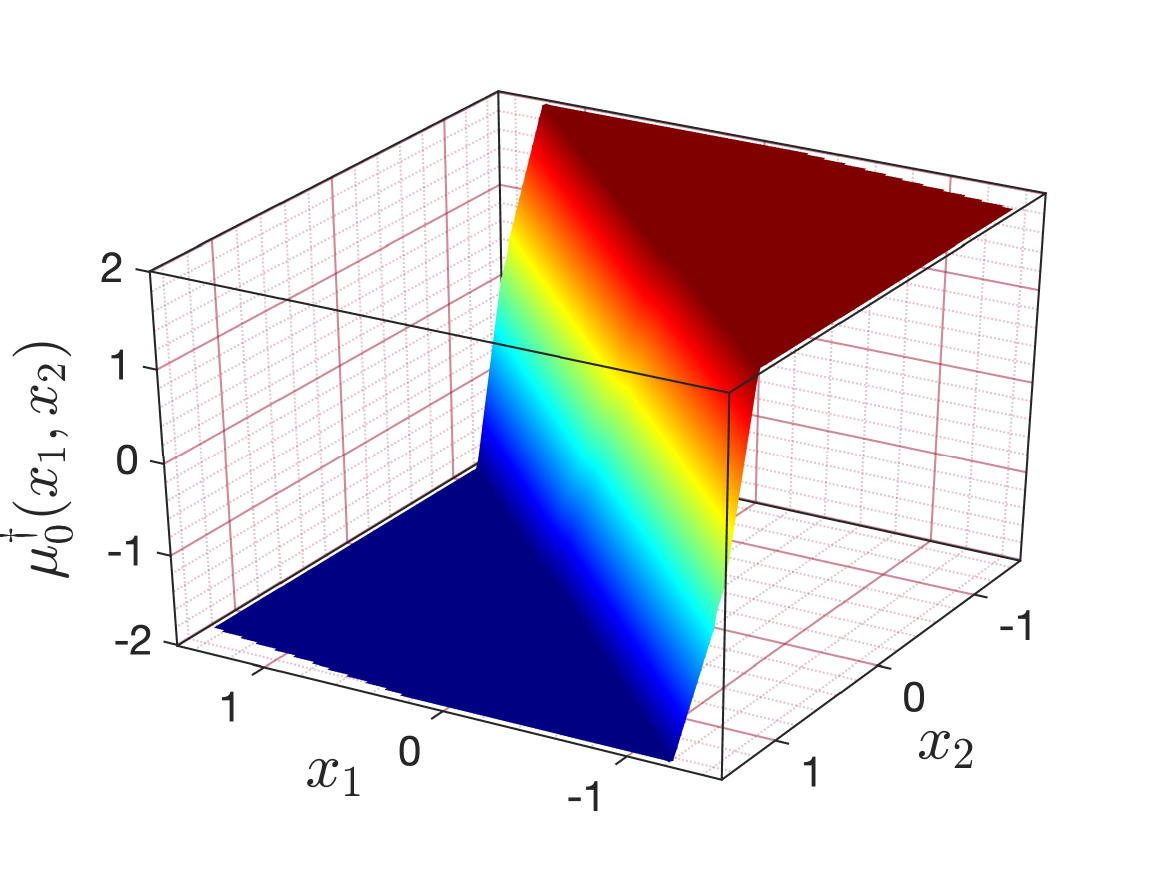}
  \end{subfigure}
  \begin{subfigure}[b]{0.49\linewidth}
    \includegraphics[width=6cm,height=5cm]{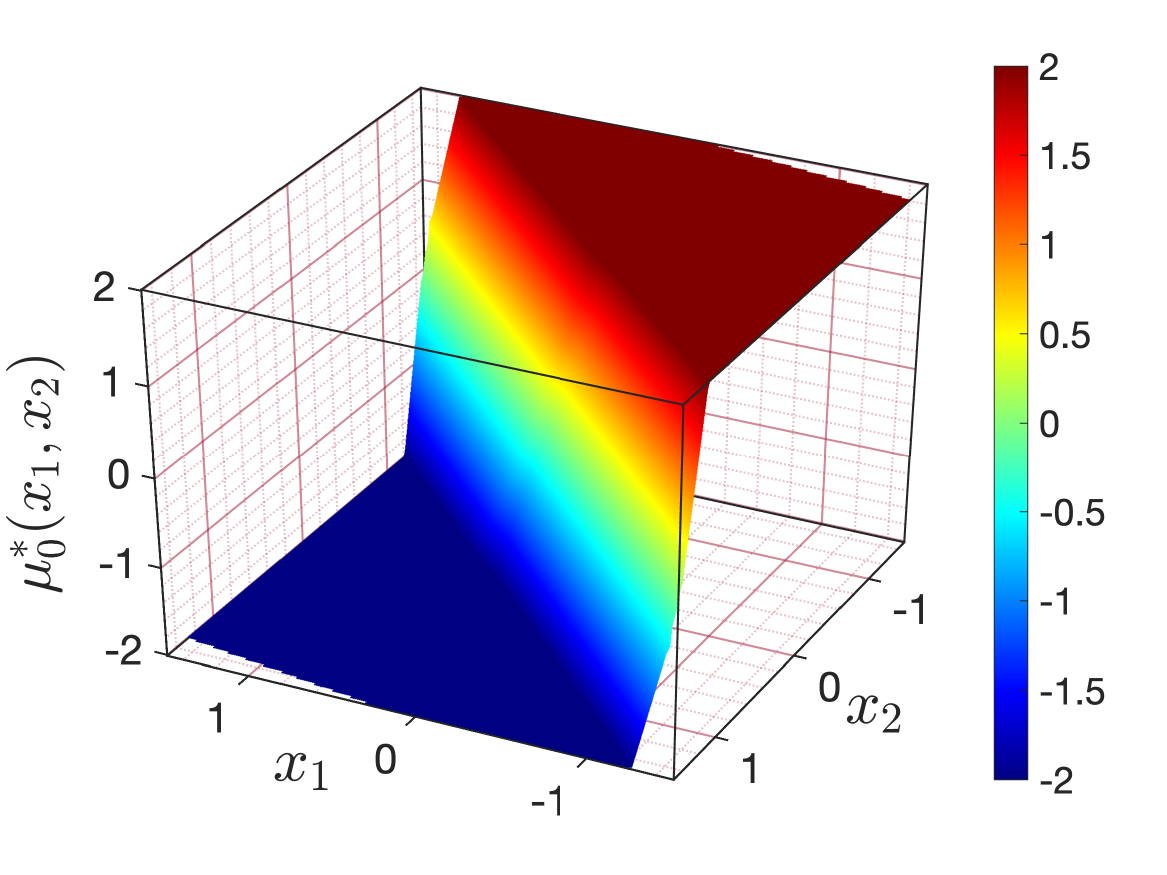}
  \end{subfigure}
\caption{The left-hand and the right-hand subfigures depicts the learned policy \(\mutrunc(\cdot)\) computed via Algorithm \ref{alg:extension_algo} and the policy \(\mu^*_0(\cdot)\), respectively.}\label{fig:uopt_uapp_rmpc_ex_1}
\end{figure}

\begin{figure}[h]
  \begin{subfigure}[b]{0.49\linewidth}
    \includegraphics[width=6cm,height=5cm]{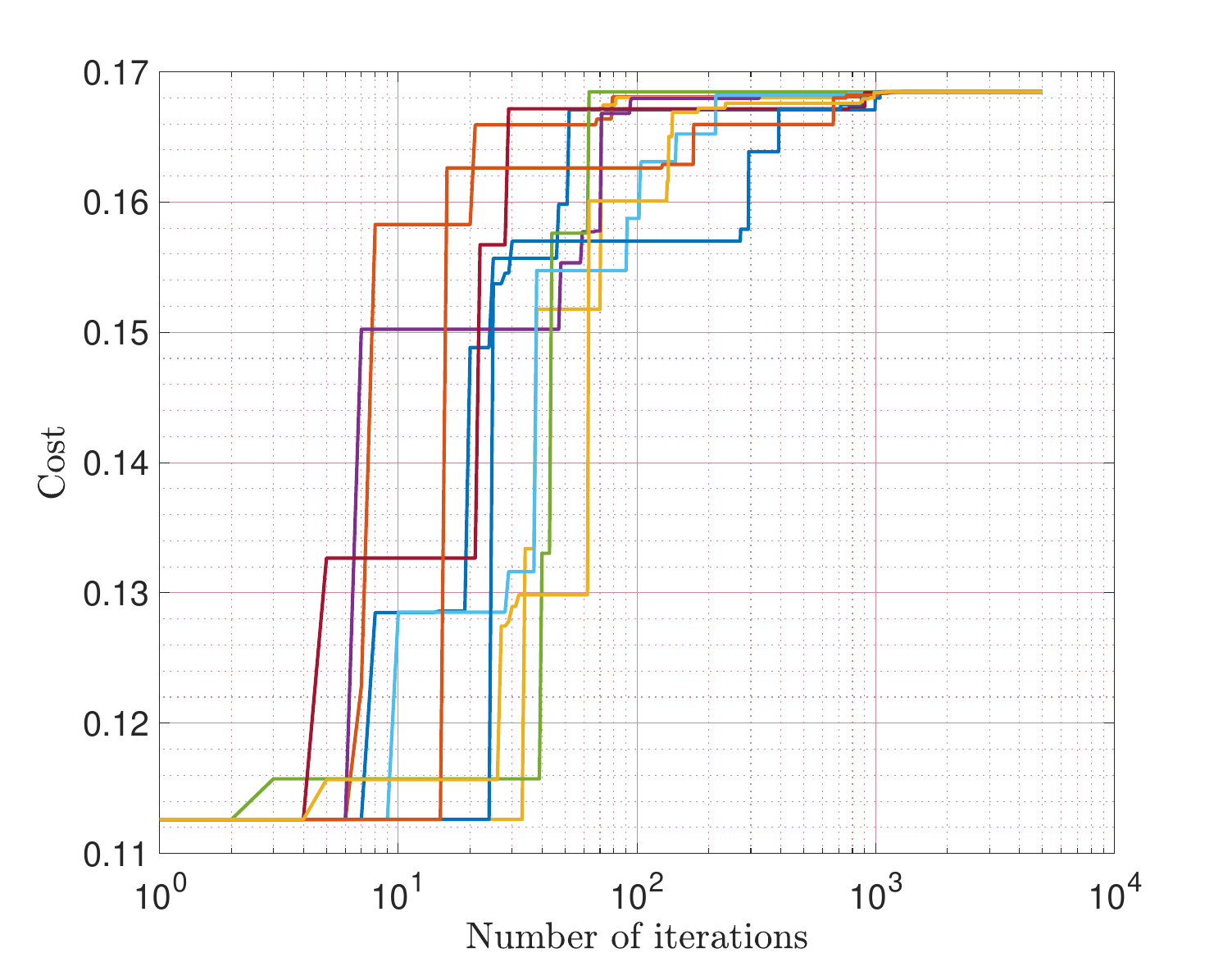}
  \end{subfigure}
  \begin{subfigure}[b]{0.49\linewidth}
    \includegraphics[width=6cm,height=5cm]{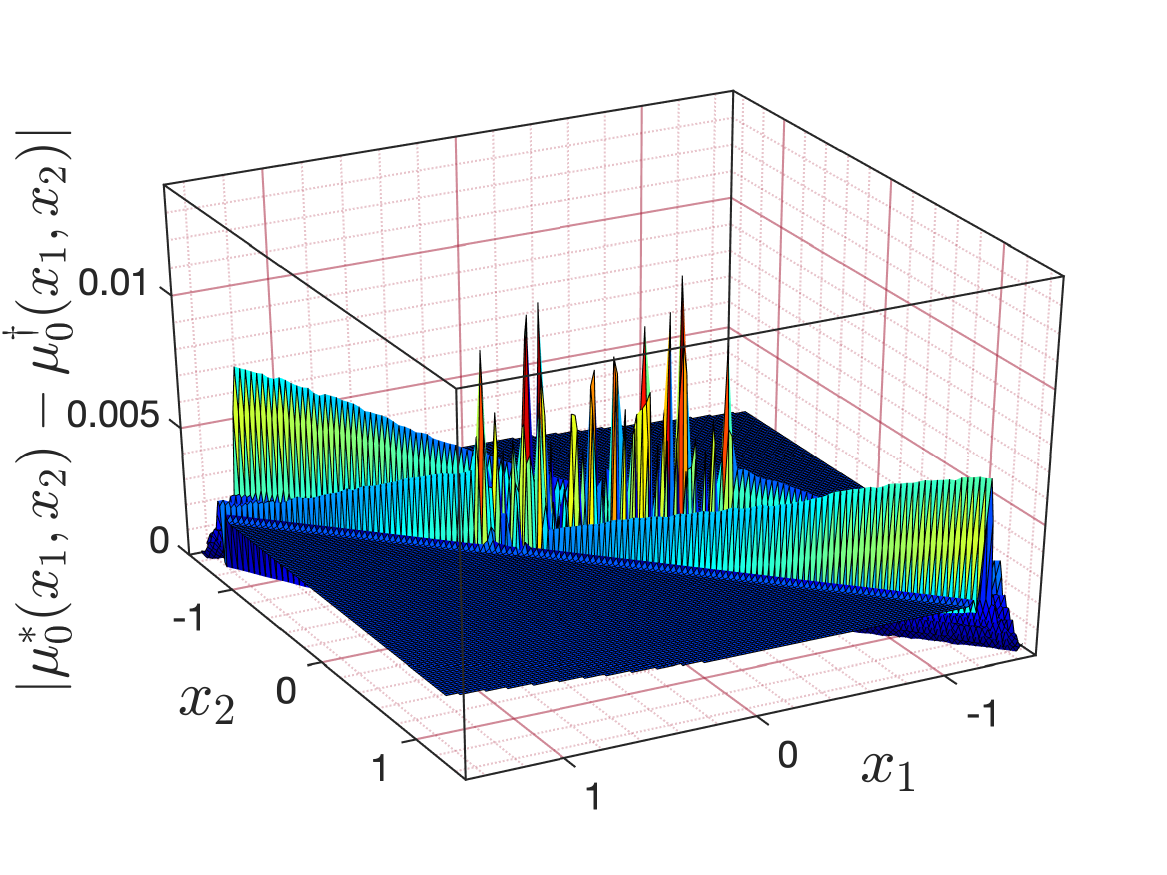}
  \end{subfigure}
\caption{The left-hand subfigure shows that the optimal value consistently converges to the same value across multiple runs of simulated annealing when Algorithm \ref{alg:msap} is employed. The right-hand subfigure depicts the error surface corresponding to the feedback policies in Figure \ref{fig:uopt_uapp_rmpc_ex_1}; it can be seen that the prespecified error margin \(\eps=0.03\) is respected}\label{fig:SipConvQError}
\end{figure}


\noindent\textbf{Learning via Algorithm \ref{alg: nn_algo}:} We move to the implementation details of the learning Algorithm \ref{alg: nn_algo} using the ReLU networks. We fixed \(\eps=0.1\) and generated \(\texttt{(state, action)}\) data by solving the SIP via Algorithm \ref{alg:msap}, independently, in parallel, at each grid point of \(\lcrc{-1.5}{1.5} \times \lcrc{-1.5}{1.5}\) with a step size of \(0.01\). The \(\texttt{(state, action)}\) data was used to train a ReLU neural network, excluding sampled initial states outside the feasible region. The final training dataset consisted of 87948 data points. Points for the testing set were sampled randomly from the feasible region amounting to a total of 10000 testing points. The neural network was trained using PyTorch on the NVIDIA Tesla P100 GPU. The input and output dimensions of the neural network were \(d = 2 \) and \(d' = 1\). For illustration, the width of the ReLU network \(W\) was fixed at 64, and using the estimates provided in \ref{thrm:main:estimate_L}, we obtained the value for the depth \( L = \ceil{22\widetilde{L} + 22} = 8721411\), where \(\widetilde{L}\) is as given in \eqref{eq:error_estimate_L}. Training such a large network is computationally challenging, but a network with a much lower depth of \(L = 4\) allows us to achieve our target error margin (verified empirically below). With these new parameters, our ReLU-network approximation \(\mutrunc(\cdot)\) (following our previous notation in \eqref{eq:nn_approx_policy}) is: 
\begin{equation}
     \mutrunc(\cdot) \Let \Upsilon^{64,4}_{\upopt(\cdot)}(\relu(\cdot);2,1).
\end{equation}
The training was done using supervised learning on the training data generated for 30000 epochs over 8 hours with a Stochastic Gradient Descent optimizer.  The learning rate was \(5\times 10^{-4}\) for all the epochs.
A typical evolution of the training loss versus the epochs is depicted in Figure \ref{fig:ValidationLoss}. Post-training, we generated the approximate feedback map using the neural network, and the result is shown in the left-hand subfigure in Figure \ref{fig:nn_approx_plus_error}. The approximation error is shown in the right-hand subfigure in Figure \ref{fig:nn_approx_plus_error}, and we note that it strictly satisfies the prespecified error-margin \(\eps.\)
\begin{figure}[h]
  \begin{subfigure}[b]{0.49\linewidth}
    \includegraphics[width=6cm,height=5cm]{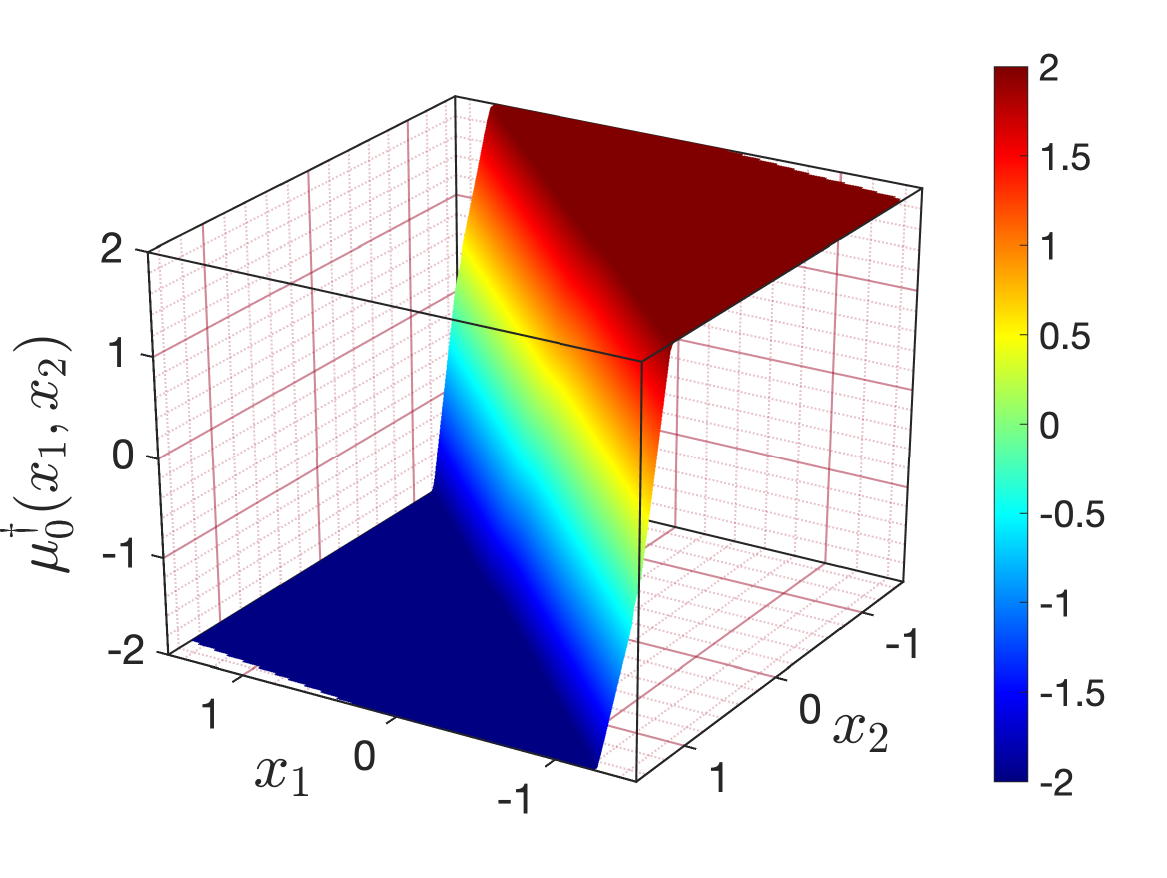}
  \end{subfigure}
  \begin{subfigure}[b]{0.49\linewidth}
    \includegraphics[width=6cm,height=5cm]{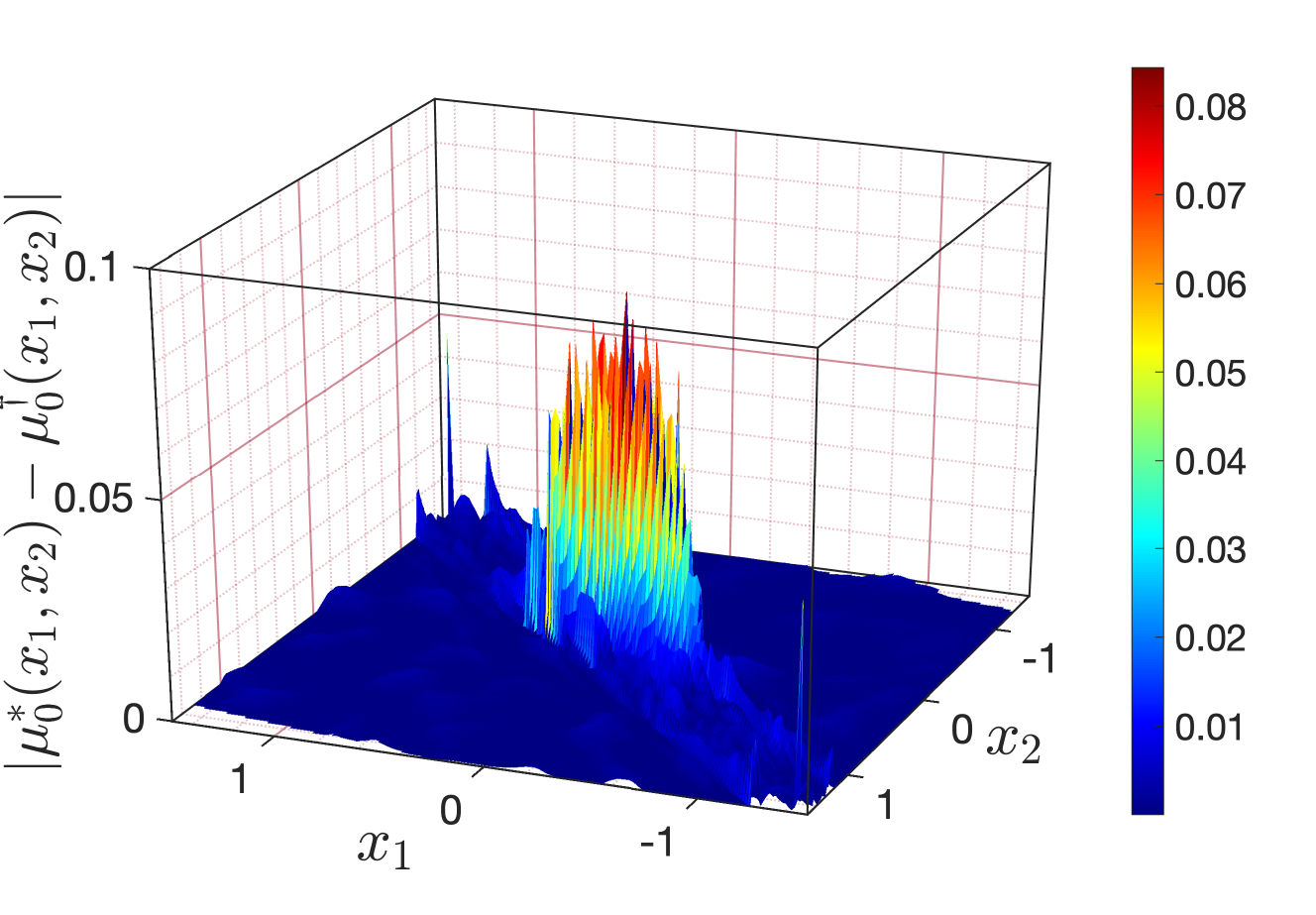}
  \end{subfigure}
\caption{The approximate policy \(\mutrunc(\cdot)\) (the left-hand subfigure) obtained via Algorithm \ref{alg: nn_algo} and the corresponding error surface (the right-hand subfigure). }\label{fig:nn_approx_plus_error}
\end{figure}


\begin{table}[ht]
\centering
\caption{Error margin and computation-time statistics. The algorithm in \cite{ref:Gao_minmax_TAC} is based on approximate multiparametric programming and there is no feature to control uniform approximation error. \textbf{CT.M.} \(\bigl(\mutrunc(\xz)\bigr)\)  represents the per-point computation time to evaluate the learned feedback map in MATLAB.}
\vspace{1mm}
\begin{tblr}{l c c}
\hline[2pt]
\SetRow{azure9}
Method  & Uniform error & \textbf{CT.M.} \(\bigl(\mutrunc(\xz)\bigr)\) (in mil. sec.) \\
\hline
\cite{ref:Gao_minmax_TAC} & NA & 15    \\
Algorithm \ref{alg:extension_algo}  & \(\eps=0.03\) & 0.2   \\
Algorithm \ref{alg: nn_algo}  &  \(\eps=0.1\) & 1.6 \\
RHC/MPC  & NA & 30  \\
\hline[2pt]
\end{tblr}
\label{tab:metadata_ex_1}
\end{table}

Table \ref{tab:metadata_ex_1} collects the computation-time data for the explicit feedback law \(\mutrunc(\cdot)\) when Algorithm \ref{alg:extension_algo}, Algorithm \ref{alg: nn_algo}, and the algorithm established in \cite{ref:Gao_minmax_TAC} were employed, respectively. Both Algorithm \ref{alg:extension_algo} and Algorithm \ref{alg: nn_algo} were computationally faster compared to the algorithm in \cite{ref:Gao_minmax_TAC}. For a fair comparison with the algorithm in \cite{ref:Gao_minmax_TAC}, we employed MATLAB to check the online evaluation time of the policies \(\mutrunc(\cdot)\) obtained via several algorithms listed in Table \ref{tab:metadata_ex_1}. However, we observed that, with Julia 10.5, the per-point computation time improved to \(13.5\) microseconds for Algorithm \ref{alg:extension_algo} and \(0.01\) microseconds for Algorithm \ref{alg: nn_algo}. This demonstrates the advantage of data-driven approximation instead of multiparametric programming-based explicit methods.
\begin{table}[ht]
\centering
\caption{Uniform and MSE error with varying depth and width of the network.}
\vspace{1mm}
\begin{tblr}{l c c c}
\hline[2pt]
\SetRow{azure9}
 Depth & Width & MSE error & Uniform Error \\
\hline
4 & 64 & \(9.620 \times 10^{-6}\) & \(0.09\) \\ 
4 & 70 & \(8.753 \times 10^{-6}\) & \(0.09\) \\
4 & 75 & \(1.020 \times 10^{-5}\) & \(0.09\) \\
5 & 64 & \(9.921 \times 10^{-5}\) & \(0.09\) \\
6 & 64 & \(1.109 \times 10^{-6}\) & \(0.09\) \\
4 & 60 & \(1.228 \times 10^{-6}\) & \(0.10\) \\
\hline[2pt]
\end{tblr}
\label{tab:nn:ex1:unif:errors}
\end{table}
As discussed earlier, due to the large network dimension estimates, we trained smaller models and verified that they satisfy the prescribed bounds; see Table \ref{tab:nn:ex1:unif:errors}. We also performed a validation process. 
\subsubsection*{Validation}
For the post-hoc validation procedure we adopt the technique established in \cite{hertneck_RMPC}. Let \(\delta_h\) be the confidence parameter and \(\tilde{\mu}\) be the empirical risk as defined in \cite[\S IV, B]{hertneck_RMPC}, then it was shown that with a confidence \(1-\delta_h\) the probability that the approximation error between the policies \(\mutrunc(\cdot)\) and \(\upopt(\cdot)\) is below the prespecified tolerance \(\eps>0\), along a trajectory with a random initial condition from \(\fset\) is larger than \(\tilde{\mu}-\eps_h\).\footnote{See \cite[Lemma 7, Equation 17]{hertneck_RMPC} for more details on several quantities related to the probabilistic validation procedure.} The validation procedure is adapted to include external disturbance sequences \(W\). For each initial point, multiple such disturbance sequences are rolled out and a constraint check is performed on all such sequences. We choose a \(\delta_h = 0.01\) and \(\mu_{\mathrm{crit}} = 0.98\). For validation rollouts, we sample initial states uniformly from the feasible set \(\fset\) with the number of initial states \(p = 50000\). The validation terminates with a \(\tilde{\mu} = 0.9895\). Hence \(\tilde{\mu} - \eps_h > \mu_{\mathrm{crit}}\) for the chosen parameters, validating the guarantees up to a confidence.
\begin{figure}
\includegraphics[scale=0.4]{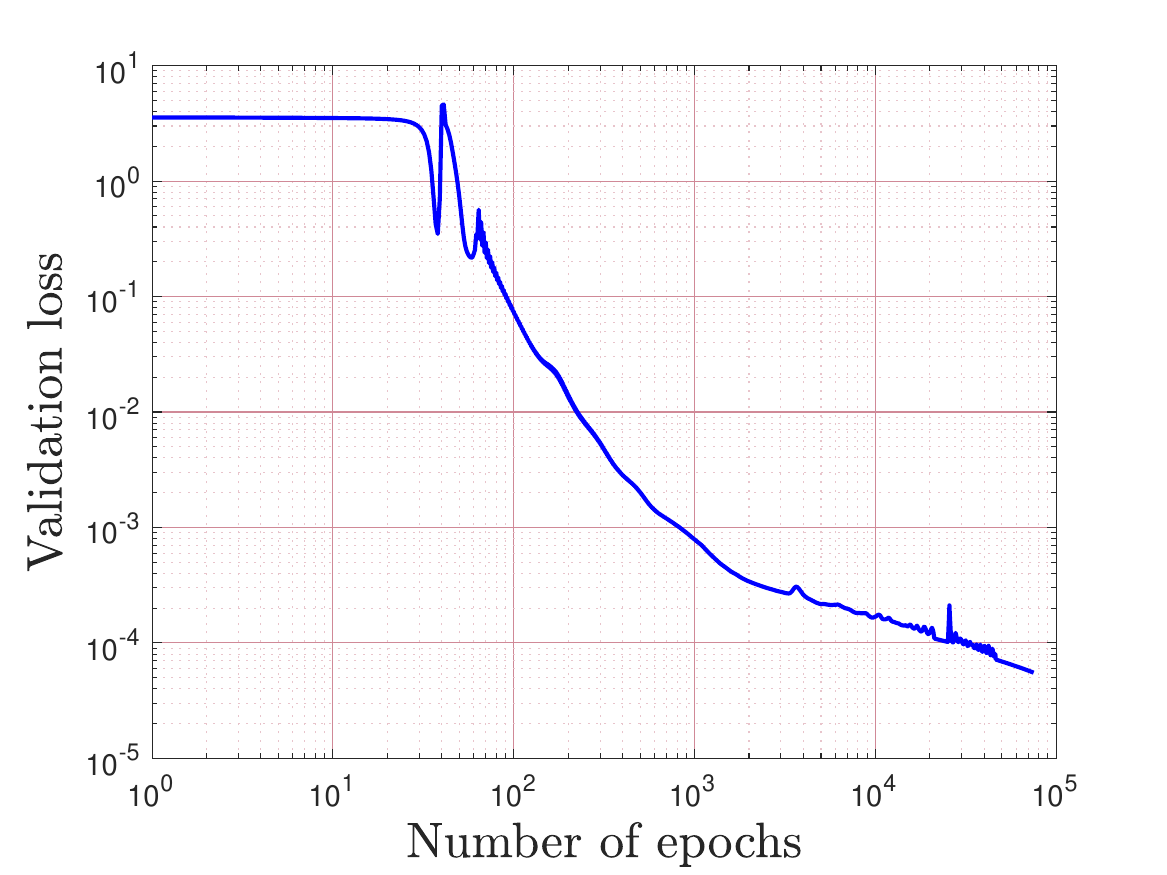}
\caption{Validation loss of the NN during training.}\label{fig:ValidationLoss}
\end{figure}
\end{example}

\begin{example}\label{exmp:lmpc_0}
Consider the fourth-order system
\begin{align}\label{eq:fourth_order_dyn_II}
    x_{t+1} = Ax_t + B u_t + G w_t,
\end{align}
with the state and the actuation matrices:
\[A \Let \begin{pmatrix} 0.40 & 0.37 & 0.29 & -0.72 \\ -0.21 & 0.64 & -0.67 & -0.04 \\ 0.83 & 0.01 & -0.28 & 0.38 \\ -0.07 & 0.60 & 0.55 & 0.49 \end{pmatrix},\,
B \Let \begin{pmatrix} 1.61 \\ 0.40 \\ -1.45 \\ -0.67 \end{pmatrix},\]
and \(G \Let (1\;1\;1\;1)^{\top}\). Fix \(N=8\) and consider the OCP
\begin{equation}
	\label{eq:RMPC_fourth_order_II}
	\begin{aligned}
		&  \min_{\pi(\cdot)} \max_{\winopt} && \hspace{-3mm}\sum_{t=0}^{\horizon-1} \inprod{\st_t}{Q\st_t} + \inprod{\cont_t}{R\cont_t} \\
		& \sbjto && \hspace{-3mm}\begin{cases}
			 \eqref{eq:fourth_order_dyn_II}, \st_0=\xz,\,\dummyx_t \in \Mbb,\text{ and }\dummyu_t \in \Ubb\\
  \text{for all }(t,\winopt)\in [0;\horizon-1]\times \Wbb^{\horizon},
		\end{cases}
	\end{aligned}
\end{equation}
where \(\Mbb \Let \aset[]{x \in \Rbb^4 \suchthat \|\dummyx\|_{\infty}\le 5}\), \(\Ubb \Let \aset[]{\dummyu \in \Rbb \suchthat |\dummyu|\le 0.2}\), and \(\Wbb \in \Rbb \Let \aset[]{\dummyw \suchthat \|\dummyw\|_{\infty} \le 0.01}\), with \(Q \Let I_{4 \times 4}\) is a \(4 \times 4\)-identity matrix, \(R \Let 0.2\). 

\noindent\textbf{Approximation via the Algorithm \ref{alg: nn_algo}:} 
Fixing \(\eps = 0.1\), we solved the ensuing SIP arising from the problem \eqref{eq:RMPC_fourth_order_II} at \(10^6\) points which were independently and uniformly randomly sampled from the state space; the computer/server/solver specifications are identical to those in Example \ref{exmp:rmpc_1}. Excluding the infeasible sample points, the final number of data points used for training were \(735105\). The input and output dimensions of the neural network were \(d = 4 \) and \(d' = 1\). The width of the ReLU network \(W\) was fixed at 256, and using the estimates provided in \ref{thrm:main:estimate_L}, we obtained the value for the depth \( L = \ceil{22\widetilde{L} + 24} = 1735359\), where \(\widetilde{L}\) is as given in \ref{eq:L_estimate}. As in Example \ref{exmp:rmpc_1},  training a network with much lower depth of \(L=6\) allowed us to achieve our error margin. With these parameters, our ReLU-network approximation \(\mutrunc(\cdot)\) (following our previous notation in \eqref{eq:nn_approx_policy}) is: 
\begin{equation}
     \mutrunc(\cdot) \Let \Upsilon^{256,6}_{\upopt(\cdot)}(\relu(\cdot);4,1).
\end{equation}
The training was done as in Example \ref{exmp:rmpc_1} using supervised learning on the training data generated for 4000 epochs with Stochastic Gradient Descent optimizer with a learning rate of \(5\times 10^{-3}\) for around 12 hours of training time on the same setup. Post training validation uniform norm error of the approximated policy over the actual feedback was \(0.0856\), i.e., 
\(
   \|\mutrunc(\cdot) - \upopt(\cdot)\|_{\mathrm{u}} \approx 0.0856.
\)
Table \ref{tab:metadata_ex_2} collects the computation-time and storage requirement data concerning the explicit feedback law \(\mutrunc(\cdot)\) when Algorithm \ref{alg: nn_algo} and the Algorithm established in \cite{ref:Gao_minmax_TAC} (via YALMIP) was employed respectively. Same as Example \ref{exmp:rmpc_1}, we observed that, with Julia 10.5, the per-point computation time improved to \(3\) microseconds for Algorithm \ref{alg: nn_algo}.
\begin{table}[t!]
\centering
\caption{Error margin and computation-time statistics where \textbf{CT.M.} \(\bigl(\mutrunc(\xz)\bigr)\) (in mil. sec.) represents the per-point computation time of the learned feedback map in MATLAB.}
\vspace{1mm}
\begin{tblr}{l c c}
\hline[2pt]
\SetRow{azure9}
Method  & Uniform error & \textbf{CT.M.} \(\bigl(\mutrunc(\xz)\bigr)\) (in mil. sec.) \\
\hline
\cite{ref:Gao_minmax_TAC} & NA & 43 \\
Algorithm \ref{alg: nn_algo} &  \(\eps=0.1\) & 2.4  \\
RHC/MPC & NA & 270 \\
\hline[2pt]
\end{tblr}
\label{tab:metadata_ex_2}
\end{table}
\begin{figure*}[h]
  \begin{subfigure}[b]{0.49\linewidth}
    \includegraphics[width=6cm,height=5cm]{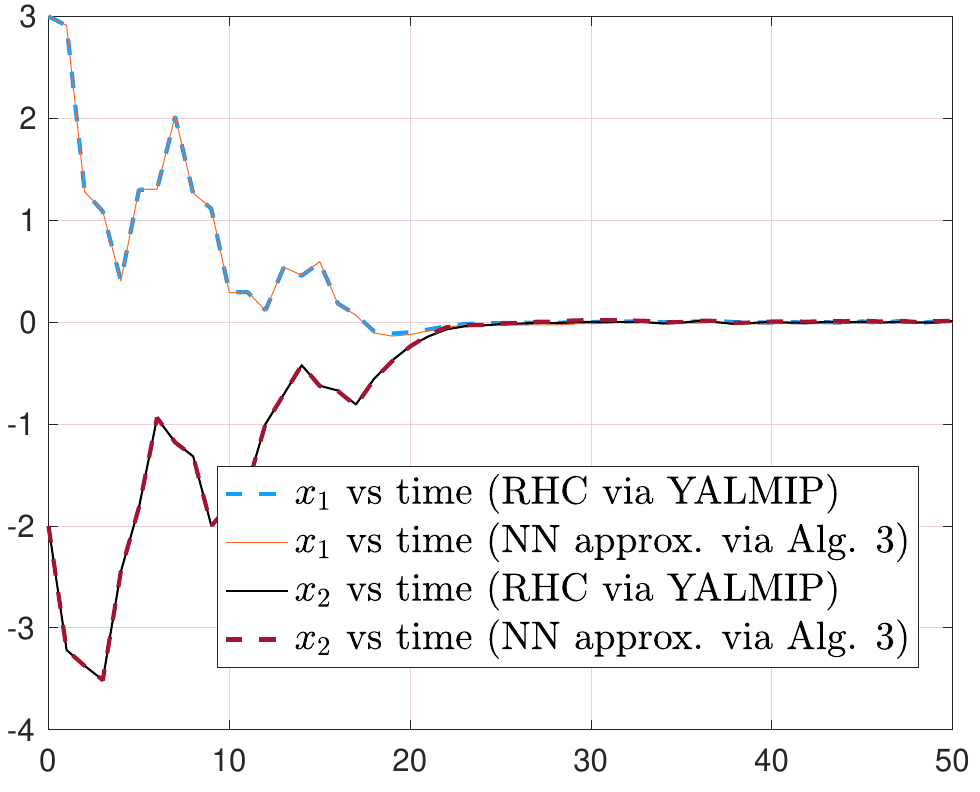}
  \end{subfigure}
  \begin{subfigure}[b]{0.49\linewidth}
    \includegraphics[width=6cm,height=5cm]{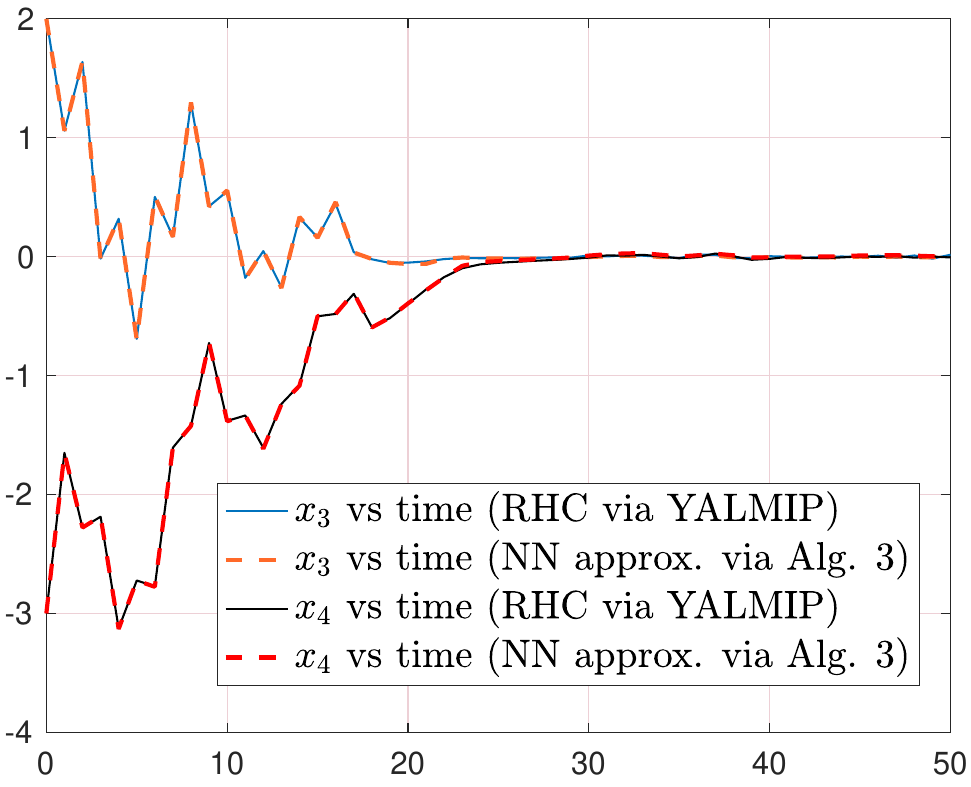}
  \end{subfigure}
\caption{Evolution of the state trajectories starting from \(\xz = (3,-2,2,-3)^{\top}\), generated using the explicit control law and the MPC/RHC control law; see Figure \ref{fig:control:laws}. }\label{fig:3rd_order_u_x1}
\end{figure*}
\begin{figure}
    \centering
    \includegraphics[scale=0.4]{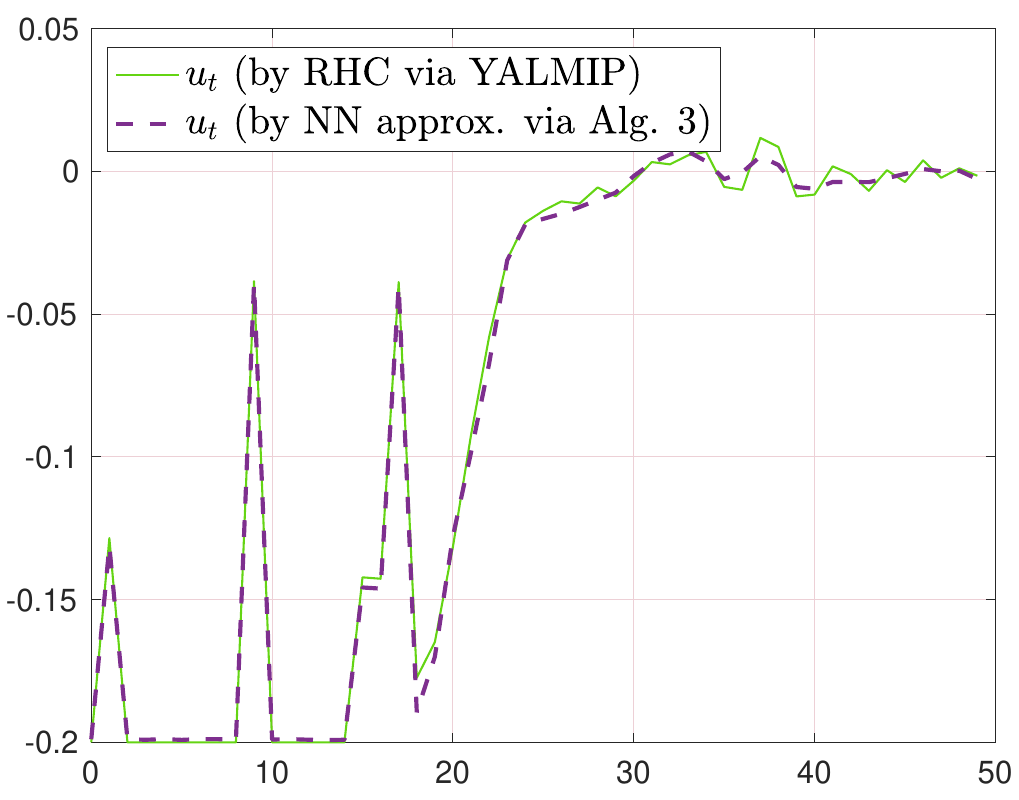}
    \caption{The explicit control law and the RHC law.}
    \label{fig:control:laws}
\end{figure}
Figure \ref{fig:3rd_order_u_x1} and \ref{fig:control:laws} depict the performance of the explicit control law generated by employing Algorithm \ref{alg: nn_algo} the corresponding closed-loop state trajectories. It can be seen that the states and the control satisfy the constraints and their trajectories almost overlap with the online receding horizon control and state trajectories with an improved computation time.

\subsubsection*{Validation}
To adopt the same post-hoc validation technique as before. We choose a \(\delta_h = 0.01\) and \(\mu_{\mathrm{crit}} = 0.99\). For validation rollouts, we sample initial states uniformly from the feasible set \(\fset\) with the number of initial states \(p = 50000\). The procedure terminates with a \(\tilde{\mu} = 0.9935\), satisfying \(\tilde{\mu} - \eps_h > \mu_{\mathrm{crit}}\) for the chosen parameters, validating the guarantees for the approximate explicit feedback law up to confidence.
\end{example}

\section{Conclusion}
In this article, we presented two feedback policy learning algorithms with uniform approximation guarantees --- one based on quasi-interpolation and the other on neural networks --- for fast computation of control policies in MPC problems. Our approach generates accurate \((\texttt{state, action})\) data by solving minmax problems exactly, representing a significant departure from traditional synthesis techniques based on multiparametric optimization or kernel methods. Several promising avenues for future work exist; for instance, one could explore the use of physics-informed neural networks or the recently introduced Kolmogorov-Arnold networks for learning approximate policies. Moreover, hardware implementations present an exciting and natural next step that we plan to investigate in subsequent studies.

\appendix

\section{Preliminaries}\label{sec:prelim}
\subsection{Stability}\label{subsec:prelim_stab}

We provide a very brief background and a definition of ISS-like stability \cite{ref:JiaWang-2001} which is needed in our study, in particular, to prove closed-loop robust stability results in Theorem \ref{thrm:stability_main_result} and Theorem \ref{thrm:stability_main_result_nn}. Recall the noisy dynamical system defined in \eqref{eq:approx-ready robust MPC}, and for any continuous admissible feedback \(x \mapsto \kappa(x) \in \Ubb\), define the vector field \(F(x_t,w_t,v_t) \Let A x_t+ B \kappa(x_t)+Bv_t+w_t\) and consider the discrete-time system 
\begin{align}\label{eq:gen:nl:sys}
    x_{t+1} = F (x_t,w_t,v_t), \quad x_0=\xz,\quad t \in \N.
\end{align}
Denote the set of sequences \(S_w\) and \(S_v\) whose elements are denoted by \(\textbf{w}\Let (w_t)_{t \in \Nz}\) and \(\textbf{v} \Let (v_t)_{t \in \Nz}\), respectively, such that \(w_t \in \Wbb\) and \(v_t \in \Vbb\) for all \(t \in \Nz\). We denote the solution of \eqref{eq:gen:nl:sys} starting from \(\xz \in \Rbb^d\), under the admissible sequence of uncertainties \(\textbf{w} \in S_w\) and \(\textbf{v} \in S_v\) by \(x(t;\xz,\textbf{w},\textbf{v})\). Also define \(W_{\text{sup}}\Let\sup_{(w,v)\in \Wbb \times \Vbb}\norm{(w,v)}\).

\begin{definition}\label{def:RobPosInv}
Consider the the system \eqref{eq:gen:nl:sys} and its corresponding notations. Let \(\Xi\) be a compact subset of \(\Rbb^d\) containing the origin \(0\in \Rbb^d\) in its interior. The set \(\Xi\) is \emph{robust positively invariant} for \eqref{eq:gen:nl:sys} if whenever \(x_0 = \xz \in \Xi\), then \(x(t;\xz,\textbf{w},\textbf{v}) \in \Xi\) for all \(t \in \N\).
\end{definition}

\begin{definition}\cite{ref:JiaWang-2001}\label{def:isps_def}
Given a compact subset \(\Xi\) of \(\Rbb^d\) containing the origin in its interior, the system \eqref{eq:gen:nl:sys} with admissible disturbance sequences \((\textbf{w},\textbf{v}) \in S_w\times S_v\) is \emph{ISS-like stable} in \(\Xi\) with respect to \(\bigl(w_t,v_t\bigr)_{t\in \Nz}\) if \(\Xi\) is robust positively invariant for \eqref{eq:gen:nl:sys}, and if there exists \(\beta(\cdot) \in \mathcal{KL}\) and \(\gamma(\cdot) \in \mathcal{K}\) such that for all  \(\xz \in \Xi\) and \(t \in \N\)
    \begin{align}
        \norm{\st(t;\xz,\textbf{w},\textbf{v})} \le \beta\bigl(\norm{\xz},t\bigr)+ \gamma \bigl(W_{\text{sup}}\bigr).
    \end{align}
\end{definition}
\subsection{A brief review of the MSA algorithm}\label{app:msa_overview}
Let \(n,k\in \N\) and consider the convex semi-infinite program (CSIP)
\begin{align}\label{eq:CSIP_original}
\inf_{\sipdvar \in \sipdvarset}  \left\{\sipcost(\sipdvar)  \;\middle\vert\;  
    \begin{array}{@{}l@{}}
        \sipcstr(\sipdvar,\sipsemvar) \le 0 \quad \text{for all}\, \sipsemvar\in \sipsemvarset
        \end{array}
        \right\}
\end{align}
where \(\sipdvarset\subset \Rbb^{n}\) is a closed and convex set such that \(\intr \sipdvarset \neq \emptyset\). The cost \(\sipcost(\cdot)\) is convex and continuous, and the constraint map \(\sipcstr(\cdot,\cdot)\) is jointly continuous, and convex in each \(\sipdvar\) for each fixed \(\sipsemvar\); \(\sipsemvarset\subset \Rbb^k\) is compact. By assumption, the interior of the feasible set \(\mathcal{H} \Let \aset[\big]{\sipdvar \in \sipsemvarset \mid \sipcstr(\sipdvar,\sipsemvar) \le 0\,\text{for all}\,\sipsemvar\in \sipsemvarset}\) is nonempty. 

A recent result presented in \cite[Theorem 1]{ref:DasAraCheCha-22} offers an efficient method for solving the CSIP \eqref{eq:CSIP_original} with minimal memory usage and serves as the foundation for the data generation process for the learning schemes in this article. 
\begin{theorem}\cite[Theorem 1]{ref:DasAraCheCha-22}\label{thrm:DACC_thrm1}
Consider the CSIP \eqref{eq:CSIP_original} along with its associated data and define \begin{align}
    \mathcal{G}(\sipsemvar_1,\ldots,\sipsemvar_{n}) \Let \inf_{\sipdvar \in \sipdvarset} \aset[\Big]{\sipcost(\sipdvar) \suchthat \sipcstr(\sipdvar,\sipsemvar_i) \le 0 \quad \text{for}\,i=1,\ldots,n}.\nn
\end{align}
Consider the global maximization problem 
\begin{align}\label{semi-inf-G}
\sup_{(\sipsemvar_1,\ldots,\sipsemvar_{n}) \in \sipsemvarset^{n}} \mathcal{G}(\sipsemvar_1,\ldots,\sipsemvar_{n}).
\end{align}
If \((\sipsemvar_1\as,\ldots,\sipsemvar_{n}\as) \in \sipsemvarset^n\) is an optimizer concerning the maximization problem \eqref{semi-inf-G}, then \(G(\sipsemvar_1\as,\ldots,\sipsemvar_{n}\as)=y\as\), where \(y\as\) is the optimal value of \eqref{eq:CSIP_original}.
\end{theorem}

\subsection{\(\lpL[p]\)-norm estimates for ReLU networks}\label{appen:NN_results}
This section collects some of the recent results on \(\lpL[p]\)-norm estimates for ReLU networks for \(p \in \lcrc{1}{+\infty}\); we shall only employ the \(p=+\infty\) case in our work. The first result (Theorem \ref{thrm:NN_approx_result_I}) establishes an explicit rate of approximation obtained via employing ReLU networks for a given function \(f \in \mathcal{C}\bigl(\lcrc{0}{1}^d\bigr)\) in terms of the width \(W_f\) and the depth \(L_f\) of the network and the dimension \(d \in \N\).
\begin{theorem}{\cite{ref:Shen_NN_2022}}\label{thrm:NN_approx_result_I}
Let \(f \in \mathcal{C}\bigl(\lcrc{0}{1}^d \bigr)\). For any \(W \in \N\), \(L \in \N\), and \(1\le p\le +\infty\), there is a function \(\widehat{f}(\cdot)\) implemented via a ReLU network with width \(C_1 \max \aset[\big]{d \lfloor W ^{1/d}\rfloor,W+2}\) and depth \(11L+C_2\) such that 
\begin{align}\label{NN_estimate_1_Cont}
    &\bigl\|f(\cdot) - \widehat{f}(\cdot)\bigr\|_{\lpL[p](\lcrc{0}{1}^d)} \le 131 \sqrt{d}  \omega_f  \Bigg(\biggl(W^2L^2 \log_3(W+2) \biggr)^{-1/d}\Bigg),
\end{align}
where \(\omega_f(r_0) \Let \sup_{x,y \in \lcrc{0}{1}^d}\aset[\big]{|f(x)-f(y)|\suchthat \|x-y\|_2 \le r_0}\) for any \(r_0 \ge 0\) is the modulus of continuity for \(f(\cdot)\). Moreover, \(C_1 = 16\) when \(p \in \lcro{1}{+\infty}\) and \(3^{d+3}\) when \(p =  +\infty\); \(C_2 = 18\) when \(p \in \lcro{1}{+\infty}\) and \(18+2d\) when \(p =  +\infty.\)
\end{theorem}
We need similar estimates for Lipschitz functions on \(\lcrc{0}{1}^d\). Observe that for Lipschitz functions \(\omega_f(t) = \lambda_f t\), where \(\lambda_f\) is the Lipschitz rank. Then the following \emph{uniform error estimate} follows immediately from Theorem \ref{thrm:NN_approx_result_I}
\begin{corollary}\label{corr:NN_approx_result_II}
Let \(f(\cdot)\) be Lipschitz with rank \(L_0\) defined on \(\lcrc{0}{1}^d\). For any \(W \in \N\) and \(L \in \N\) with \(W \ge 3^{d+4}d\) and \(L \ge 29+2d\), there is a function \(\widehat{f}(\cdot)\) implemented via a ReLU network with width \(W\) and depth \(L\) such that 
\begin{align}\label{eq:NN_estimate_2_Lip}
    \unifnorm{f(\cdot) - \widehat{f}(\cdot)}\le& 131\sqrt{d} L_0 \Bigg{(}\biggl(\tfrac{W}{3^{d+5}d}\biggr)^2 \biggl( \tfrac{L-18-2d}{22}\biggr)^2  \log_3 \biggl( \tfrac{W}{3^{d+5}d}+2 \biggr) \Bigg{)}^{-1/d}.
\end{align}
\end{corollary}

\section{Proofs of Lemma \ref{lem:Aux_lem} and Theorem \ref{thrm:value_func_equality}}\label{appen:msa:proofs}
\begin{prooflemm}
Consider the SIP \eqref{eq:approx ready param SIP robust MPC}; for a fixed \(W \in \agradmdist^{\horizon}\), let us define the set 
\[
\feasSip_{\sipfeasbset} \Let \left\{(\theta,\eta,r) \;\middle\vert\;  
\begin{array}{@{}l@{}}
        \text{Constraints in \eqref{eq:approx ready param SIP robust MPC} hold}
        \end{array}
        \right\}.    
\]
Notice that \(\feasSip_{\sipfeasbset}\) is closed and convex for every \(\widetilde{w}_t\in \agradmdist\) for all \(t=0,\ldots,\horizon-1.\) Indeed, \((\theta,\eta,r) \mapsto \mathbb{J}_{\horizon}(\xz,\theta,\eta,W)-r\) is convex. The mapping from \((\theta,\eta,r)\) to the solution of the recursion \eqref{eq:noisy-system}, i.e., \((\theta,\eta,r) \mapsto x_{\theta,\eta}(t;\xz,W)\) and the map \((r,\theta,\eta) \mapsto (\cont_0,\ldots,\cont_{\horizon-1}) \teL U = \theta W+\eta\) is linear in \((\theta,\eta)\) and thus convex \cite[Chapter 7]{ref:Lof-03}. Thus the intersection
\(\feasSip \Let \bigcap_{W \in \agradmdist^{\horizon}} \feasSip_{\sipfeasbset}\)
is also closed and convex since it is the arbitrary intersection of closed and convex sets.  \QEDA
\end{prooflemm}

\begin{prooftheo}
Let us show that \(\gfunc(\cdot;\xz)\) is upper semicontinuous for every \(\xz \in \fset\). Consider the relaxed problem \eqref{eq: inner param SIP robust MPC} corresponding to the original SIP \eqref{eq:approx ready param SIP robust MPC}, and define the feasible set associated with the \eqref{eq: inner param SIP robust MPC} by
\begin{align}
    \mathsf{F}(\mathcal{W};\xz) \Let \hspace{-1mm} \left\{(\theta,\eta,r) \;\middle\vert\;  
    \begin{array}{@{}l@{}}
\mathbb{J}_{\horizon}\bigl(\xz,\theta^i,\eta^i,\conin{i}) - r \leqslant 0,\; \hspace{-1mm}\st_{\horizon}\in \admfinst,\, \st(0)=\xz,\\ \st^i_{\theta,\eta}(t;\xz,\conin{i}) \in \Mbb,\, \cont^i_t+v^i_t \in \Ubb,\,\conin{i} \in \agradmdist^{ \horizon} 
        \end{array}
        \right\}\nn
\end{align}
for \(i=1,\ldots,\dvar\). The problem \eqref{eq: inner param SIP robust MPC} is now translated to
\begin{equation}\label{eq:modified_g_func}
\gfunc(\mathcal{W}; \xz)  \Let \inf_{(\theta,\eta,r)\in \sipfeasbset}\,\aset[\big]{ r \suchthat (\theta,\eta,r) \in \mathsf{F}(\mathcal{W};\xz)}.
\end{equation}
Observe that:
\begin{enumerate}[label=\textbf{(\textsc{P}\arabic*)}, leftmargin=*]
\item \label{cont_1} the map \((\theta,\eta,r) \mapsto \mathbb{J}_{\horizon}(\xz,\theta,\eta,W)-r\) is continuous in \((\theta,\eta,r)\); 
\item \label{cont_2} the map \((\theta,\eta,r) \mapsto \theta W+\eta \) is continuous in \((\theta,\eta,r)\); 
\item \label{cont_3} the map \((\theta,\eta,r) \mapsto \st_{\theta,\eta}(t;\xz,W)\) is continuous in \((r,\theta,\eta)\) for all \(t=0,\ldots,\horizon-1\).
\end{enumerate}
Fix \((\ol{\theta},\ol{\eta},\ol{r}, \ol{\mathcal{W}}) \in \sipfeasbset \times \agradmdist^{\totdvar}\), where \((\ol{\theta},\ol{\eta},\ol{r}) \in \mathsf{F}(\ol{\mathcal{W}};\xz)\). Invoking Assumption \ref{assum:slater's} and the continuity of the maps in \ref{cont_1}--\ref{cont_3}, there exists a sequence \(\bigl(\ol{\theta}_n,\ol{\eta}_n, \ol{r}_n\bigr)_{n \in \Nz} \subset F(\ol{\mathcal{W}};\xz)\) such that \(\bigl(\ol{\theta}_n,\ol{\eta}_n, \ol{r}_n\bigr)_{n \in \Nz}  \to \bigl(\ol{\theta},\ol{\eta}, \ol{r}\bigr)\) as \(n\to +\infty\). This immediately implies that the constraint qualification condition (CQ) in \cite[Definition 5.3, pp 53]{ref:param_opt_still} is satisfied, i.e., for every \(\xz \in \fset\) the following hold:
\(
(\ol{\theta}_n,\ol{\eta}_n,\ol{r}_n) \in \intr{\mathsf{F}(\ol{\mathcal{W}};\xz)}.
\)
Upper semicontinuity of \(\agradmdist^{\totdvar} \ni  \mathcal{W} \mapsto \gfunc(\mathcal{W}; \xz) \in \Rbb\) now follows readily by invoking \cite[Lemma 5.4-(b), p. 54]{ref:param_opt_still}.

The existence of an optimizer \(\mathcal{W}\as \in \agradmdist^{\totdvar}\) follows immediately from the upper semicontinuity of \(\gfunc(\cdot;\xz)\) from \ref{thrm:value_func_equality_0}, compactness of \(\agradmdist^{\dvar}\), and from the Weierstrass Theorem, thereby proving \ref{thrm:value_func_equality_1}. 

From the Assumption \ref{assum:slater's} there exists a pair \((\theta,\eta,r) \in \sipfeasbset\) such that for all \(\mathcal{W} \in \agradmdist^{\totdvar}\), the interior of the set 
 \(
 \bigcap_{(\conin{i})_{i\in [1;\dvar]} } \feasSip_{\sipfeasbset}
 \)
is nonempty. Moreover, in the decision and semi-infinite variable pair, the constraint maps are jointly continuous. Indeed the maps \((\theta,\eta,r,W) \mapsto x_{\theta,\eta}(t;\xz,W)\) and the map \((r,\theta,\eta,W) \mapsto (\cont_0,\ldots,\cont_{\horizon-1}) \teL U = \theta W+\eta\) are continuous and consequently \((\theta,\eta,r,W) \mapsto \mathbb{J}_{\horizon}(\xz,\theta,\eta,W)-r\) is continuous (recall that \(c(\cdot,\cdot)\) was jointly continuous). Then Lemma \ref{lem:Aux_lem}, the problem data \ref{eq:system-st-con-dist}--\ref{eq:system:terminal_set} and Theorem \ref{thrm:DACC_thrm1} assures that there exists an optimizer \(\mathcal{W}\as \in \agradmdist^{\totdvar}\) such that 
 \[
  \mathbb{J}\as_{\horizon}(\xz) = \gfunc(\mathcal{W}\as; \xz)\,\,\text{for every }\xz\in \fset,
 \]
which completes the proof. \QEDA
\end{prooftheo}

\section{Proofs of Theorem \ref{thrm:stability_main_result} and Theorem \ref{thrm:stability_main_result_nn}} \label{appen:ExMPC:proofs}
\begin{proof_t1}
We begin by giving a proof of the first assertion. Fix \(\eps>0\) and \(\psi(\cdot) \in \mathcal{S}(\Rbb^d)\) satisfying moment condition of order \(M\) and decay condition of order \(K>d\), with \(C_0\) as the upper bound in \eqref{eq:decay_condition}. Recall the approximation schemes \eqref{eq:approx_policy_parent_mainres} and \eqref{eq:approx_policy}. Then we have the estimate (see \eqref{eq:lipschitz_estimate_mpc_proof_mainres}) 
\begin{align}\label{eq:lipschitz_estimate_mpc_proof}
    	\unifnorm{\apprfb(\cdot) - \upopt(\cdot)} \le C_{\gamma}L_{0}h\sqrt{\Dd}+ \Delta_0(\psi,\Dd) \,\,\text{on}\,\,\Rbb^d,
    \end{align}
where \(\upopt(\cdot)\) is the extended policy \eqref{eq:extented_policy} (which is Lipschitz continuous with rank \(L_0\)), the quantity \(\Delta_0(\psi,\Dd) \Let \mathcal{E}_0(\psi,\Dd)\unifnorm{\upopt(\cdot)}\) is the saturation error, the term \(\mathcal{E}_0(\psi,\Dd)\) and \(C_{\gamma} \) are defined below \eqref{eq:lipschitz_estimate_mpc_proof_mainres}. From \cite[Chapter 2, Corollary 2.13]{ref:mazyabook} it follows that for the preassigned \(\eps>0\), we can find \(\Dd_{\mathrm{min}}>0\) such that whenever \(\Dd \ge \Dd_{\mathrm{min}}\), we have \(\mathcal{E}_0(\psi,\Dd) \le \frac{\eps}{3\unifnorm{\upopt(\cdot)}}\). Fixing such a \(\Dd \ge \Dd_{\mathrm{min}}\), we ensure that \(\Delta_0(\psi,\Dd) \le \frac{\eps}{3}.\)
Fixing 
\begin{align}
    	\label{h_estimate}
    	h = \frac{\eps}{3C_{\gamma}L_{0}\sqrt{\Dd}},
\end{align}
from \eqref{eq:lipschitz_estimate_mpc_proof} we arrive immediately arrive at \(\unifnorm{\upopt(\cdot)-\appr{\mu}_0(\cdot)} \le \frac{2\eps}{3}.\) We also have the estimate \cite[\S 2.3.2, eq. (2.57)]{ref:mazyabook}
\begin{align}\label{eq:truncated_bound_1}
    \| \mutrunc(x)-\widehat{\mu}_0(x) \| \le \mathcal{B} \biggl( \tfrac{\sqrt{\Dd}}{\rzero} \biggr)^{K-d} \unifnorm{\widehat{\mu}_0(\cdot)} 
\end{align}
for all \(x \in \Rbb^d\), where \(\mathcal{B}\) is a constant specific to \(\psi(\cdot)\). In view of \eqref{eq:truncated_bound_1} we pick
\[ \rzero \Let \sqrt{\Dd}\biggl(\frac{\eps}{3\mathcal{B} \unifnorm{\upopt(\cdot)}}\biggr)^{1/(d-K)}.\]
Then \(\|\appr{\mu}_0(x) - \mutrunc(x)\| \le \frac{\eps}{3}\) for all \(x \in \Rbb^d\). Now restricting the domains of \(\upopt(\cdot)\), \(\appr{\mu}_0(\cdot)\), and \(\mutrunc(\cdot)\) to \(\fset\) while retaining the same notation for all of them, we see that
\begin{align}\label{pf:final_error_estimate}
      & \bigl\|\upopt(x)-\mutrunc(x)\bigr\|   \le \bigl\|\upopt(x)-\appr{\mu}_0(x)\bigr\| + \bigl\| \appr{\mu}_0(x) - \mutrunc(x)\bigr\| \\ & \le C_{\gamma}L_0h\sqrt{\Dd}+\Delta_0(\psi,\Dd)+\mathcal{B} \biggl(\frac{\sqrt{D}}{\rzero}\biggr)^{K-d}\unifnorm{\upopt(\cdot)} \nn \le \tfrac{2\eps}{3}+\tfrac{\eps}{3}= \eps.
\end{align}
In summary, since \(\eps>0\) was preassigned and we picked \(\bigl(h,\Dd,\rzero\bigr) \in \loro{0}{+\infty}^3\) such that the estimate \eqref{pf:final_error_estimate} holds, the first assertion \ref{thrm:main:estimate} stands established. 
\par We proceed to prove the second assertion \ref{thrm:main:stability} concerning stability of the closed loop process \( (\st_t)_{t \in \Nz}\) corresponding to the system \eqref{eq:system} under the approximate feedback \(\mutrunc(\cdot)\). Under \(\mutrunc(\cdot)\), the closed loop process is given by:
\begin{align}\label{proof:closed_loop_1}
\dummyx^{+} = A\dummyx + B \mutrunc(\dummyx) + w,
\end{align}
	where \(\dummyx^+\) denotes the next state when the current state is \(\dummyx\).  Recall that the (state-dependent) approximation noise is given by \(v \Let \mutrunc(\dummyx)- \upopt(\dummyx)\). Then from \eqref{proof:closed_loop_1} we have
\begin{align}\label{proof:closed_loop_2}
\dummyx^{+} & = A \dummyx+B \mutrunc(\dummyx)-B \upopt(\dummyx)+B\upopt(\dummyx) + w = A \dummyx + B \bigl(\mutrunc(\dummyx) -\upopt(\dummyx)\bigr) + B \upopt(\dummyx) + w \nn \\& = A \dummyx + B \upopt(\dummyx) + B \dummyv + w.
\end{align}
Notice that the closed-loop dynamics \eqref{proof:closed_loop_2} is essentially the noisy system \eqref{eq:noisy-system} under the policy \(\upopt(\cdot)\), which, under Assumption \ref{assump:stability_assumptions} is stable in the sense of Definition \ref{def:isps_def} and the ensuing optimal control problem \eqref{eq:approx ready param robust MPC} is recursively feasible; see \cite[\S3]{ref:MayFal-19}. This immediately proves ISS-like stability of the system \eqref{eq:system} under the approximate feedback policy \(\mutrunc(\cdot)\) under the equivalence established in \eqref{proof:closed_loop_2} (i.e., stability of the system \eqref{proof:closed_loop_1}) in the sense of Definition \ref{def:isps_def}. \QEDA
\end{proof_t1}

\begin{proof_t2}
We first prove the error estimate \ref{thrm:main:estimate_L}, and to this end, fix \(\eps>0\). Recall that the policy \(\upopt(\cdot)\) is Lipschitz continuous and \(L_0\) is its minimal Lipschitz rank, and the ReLU-NN approximate feedback policy is given by \(\mutrunc(\cdot)\). Recalling the estimate \eqref{eq:NN_estimate_2_Lip} in Corollary \ref{corr:NN_approx_result_II}, we introduce the simplified notations \(\widetilde{L} \Let \frac{L - 18 - 2d}{22}\), where \(L\) denotes the depth of the neural network, and \(\widetilde{W} \Let \frac{W}{3^{d+5}d}\), where \(W\) represents the network's width. With these definitions, we reformulate the estimate \eqref{eq:NN_estimate_2_Lip}, replacing \(f(\cdot)\) and \(f_{\relu}(\cdot)\) with \(\upopt(\cdot)\) and \(\mutrunc(\cdot)\), respectively. Then we have the estimate
\begin{equation}
        \label{eq:nn_approx_error_proof}
        \unifnorm{ \upopt(\cdot) - \mutrunc(\cdot)} \hspace{-1mm}\le 131\sqrt{d} L_0 \left(  \widetilde{W}^2 \widetilde{L}^2 \log_3(\widetilde{W} + 2) \right)^{-\frac{1}{d}},
    \end{equation}
To ensure that the right-hand side of \eqref{eq:nn_approx_error_proof} remains bounded by \(\eps\), we must appropriately determine the parameters within \(\widetilde{L}\) and \(\widetilde{W}\). Since the dimension \(d\) is fixed, assigning a specific value \(W = \mathcal{W}_0\) to the width of the ReLU network uniquely determines \(\widetilde{W}\). We then proceed by selecting
     \begin{equation}
        \label{eq:L_estimate}
         \widetilde{L} \Let \sqrt{ \frac{\eps^{-d}}{(131\sqrt{d}L_0)^{-d}\widetilde{W}^2 \log_3(\widetilde{W} + 2)}}.
     \end{equation}
Substituting \eqref{eq:L_estimate} in \eqref{eq:nn_approx_error_proof}, we get
    \begin{align}
        \label{eq:error_estimate_L}
       & \unifnorm{ \upopt (\cdot) - \mutrunc (\cdot) } \le 131\sqrt{d} L_0  \left(  \widetilde{W}^2 \frac{\eps^{-d} \log_3(\widetilde{W} + 2)}{\bigl(131\sqrt{d}L_0\bigr)^{-d}\widetilde{W}^2 \log_3(\widetilde{W} + 2)}  \right)^{-1/d} = \eps,
\end{align}
and the value for \(L\) can then be calculated as \( L = \ceil{22\widetilde{L} + 18 + 2d}\). Thus, given a preassigned \(\eps >0\) and a fixed width \(W=\mathcal{W}_0 \in \N\), we find \(L \in \N\) such that the criterion of uniform \(\eps\)-tolerance is respected, and the assertion \ref{thrm:main:estimate_L} follows.
    
We now proceed to prove assertion \ref{thrm:main:estimate_W}. In the same spirit, given \(L=\mathcal{L}_0\) for the depth of the neural network, \(\widetilde{L} = \frac{\mathcal{L}_0-18-2d}{22}\) gets fixed. Now we set 
     \begin{equation}
        \label{eq:W_estimate}
         \widetilde{W} \Let 3^{\biggl[ \tfrac{\eps^{-d}}{(131\sqrt{d}L_0)^{-d}\widetilde{L}^2}\biggr]^{1/3}} -1.
     \end{equation}
By leveraging the bounds \( \log_3(x+1) \leq x \) and \( \log_3(x+2) \geq \log_3(x+1) \), we can reformulate \eqref{eq:nn_approx_error_proof} in an alternative form. This allows us to express the error estimate in a manner that facilitates further analysis. Specifically, we rewrite it as
\begin{align}
        \label{eq: ineq_erro_nn}
        &\unifnorm{ \upopt(\cdot)- \mutrunc(\cdot)} \le 131\sqrt{d} L_0 \left(\left(  \widetilde{W}^2 \widetilde{L}^2 \log_3(\widetilde{W} + 2) \right)^{-\frac{1}{d}}\right) \nn \\  
        &\leq 131\sqrt{d} L_0 \left(\left(  \widetilde{L}^2 \big[\log_3(\widetilde{W} + 1)\big]^3 \right)^{-\frac{1}{d}}\right),
    \end{align}
and substituting \eqref{eq:W_estimate} in \eqref{eq: ineq_erro_nn}, we get
\begin{align}
\label{eq:error_estimate_W}
& \unifnorm{ \upopt(\cdot) - \mutrunc(\cdot)} \le 131\sqrt{d} L_0  \left( \left( \frac{\eps^{-d}}{(131\sqrt{d}L_0)^{-d}\widetilde{L_0}^2}\right) \widetilde{L_0}^2 \right)^{-1/d}= \eps.
\end{align}
From \(\widetilde{W}\) we find the width \(W = \ceil{3^{d+5}\,d\,\widetilde{W}}\) of the neural network. We have shown that given \(\eps >0\) and a fixed \(\mathcal{L}_0 \in \N\), we can find a \(W \in \N\) such that the given tolerance is respected, and assertion \ref{thrm:main:estimate_W} stands established.  
    
Finally, \ref{thrm:main:stability_nn} follows verbatim from the proof in \S\ref{sec:main_tech_results} because the uniform error estimates \eqref{eq:error_estimate_L} and \eqref{eq:error_estimate_W} hold.
\QEDA
\end{proof_t2}

\bibliographystyle{amsalpha}
\bibliography{refs}

\end{document}